%% file: F-ODEs.tex
\documentclass{siamart0516}

\input{ex_shared-2}

\ifpdf
\hypersetup{
  pdftitle={\TheTitle},
  pdfauthor={\TheAuthors}
}
\fi

\begin{document}

\maketitle

\begin{abstract}
In this work, we study the long time behaviors, including asymptotic contractivity and dissipativity,  
of the solutions to several numerical methods for fractional ordinary differential equations (F-ODEs).
The existing algebraic contractivity and dissipativity 
rates of the solutions to the scalar F-ODEs are first improved.
In order to study the long time behavior of numerical solutions to 
fractional backward differential formulas (F-BDFs), 
two crucial analytical techniques are developed, with the first one 
for the discrete version of the fractional generalization of the traditional Leibniz rule,
and the other for the algebraic decay rate of the solution to a linear Volterra difference equation. 
By mens of these auxiliary tools and some natural conditions,  
the solutions to F-BDFs are shown to be contractive and dissipative, 
and also preserve the exact contractivity rate of the continuous solutions.
Two typical F-BDFs, based on the Gr\"{u}nwald-Letnikov formula and L1 method respectively,  are studied. 
For high order F-BDFs, including some second order F-BDFs and $3$-$\alpha$ order method, 
their numerical contractivity and dissipativity are also developed under some slightly stronger conditions.
Numerical experiments are presented to validate the long time qualitative characteristics
of the solutions to F-BDFs, revealing very different decay rates of the numerical solutions in terms of the  
the initial values between F-ODEs and integer ODEs and 
demonstrating the superiority of the structure-preserving numerical methods. 
\end{abstract}



\section{Introduction}
Fractional calculus has been widely applied to many areas in science and engineering.
Various fractional-order dynamical models have been proposed in applications,
and their numerical solutions have shown better consistencies with experimental data 
than those produced by the corresponding integer-order differential equations \cite{Kil, Pet, Pod}.
A typical model is the time fractional anomalous diffusion model,   
which describes a diffusion process where the mean square displacement of a particle grows slower or faster than that in the normal diffusion process. 
Anomalous diffusions were observed and confirmed in many experiments. 
Solutions to fractional anomalous diffusion equations demonstrate a very important feature, i.e., 
they decay asymptotically in an algebraic decay rate, resulting in the so-called long-tail effect.
It is highly interesting and important both mathematically and practically 
if we could have a quantitative understanding of the long-time dynamical behaviors of 
the solutions to nonlinear fractional models, especially of how the numerical solutions decay and 
if they can preserve the exact same algebraic decay rate as their continuous counterparts. 
This is a challenging topic and has basically still not been investigated in the literature, and 
will be the main motivation and focus of the current work.
Let us start with the model of our main interest.
For $0<\alpha<1$, we consider the Caputo F-ODEs:
\begin{equation} \label{eq1}
\begin{split}
~^{C}_{0}D_{t}^{\alpha}x(t)=f(t, x(t)),~~~x\in\mathbb{R}^{d},
\end{split}
\end{equation}
with initial condition $x(0)=x_{0}$, where $^{C}_{0}D_{t}^{\alpha}x(t)$ is the Caputo fractional derivative:
\begin{equation*} 
\begin{split}
^{C}_0{D}^\alpha_tx(t)=
\frac{1}{\Gamma(1-\alpha)}\int_0^t\frac{x^{(1)}(\tau)}{(t-{\tau})^{\alpha}}d\tau,~~~~t>0.
\end{split}
\end{equation*}

The stability analysis of F-ODEs has attracted a great attention and the main difficulty in the analysis 
lies in the nonlocal nature of fractional derivatives.
A fundamental stability result for linear F-ODEs, i.e., $f(t, x)=Ax$ in (\ref{eq1}), was established 
by Matignon \cite{Matignon},  where the stability region and  a concrete algebraic decay rate, namely $O(t^{-\alpha})$, 
of the solutions were derived.  Many important results and various analytical strategies for the stability of fractional linear systems 
 have been developed in succession; see the survey article \cite{LiCP}.
 
For the stability of nonlinear F-ODEs, a popular approach is to extend the
classical Lyapunov theorem to fractional systems and make use of the fractional comparison principle.
The concept of the Mittag-Leffler stability and the fractional Lyapunov second method were developed in \cite{Mittag}. 
This method relies on an appropriate Lyapunov function and the calculation of the Caputo fractional derivative
of the function. 
Under the classical Lipschitz hypothesis on function $f$, 
the stability with respect to initial values and the structural stability of F-ODEs were studied in \cite{DF}.
The stability theory of nonlinear F-ODEs is still far from maturity due to 
the coupling between the complex structure of the nonlinear function $f$ 
and the nonlocal feature of fractional derivatives.
To illustrate the motivation of the contractivity, also called the stability or asymptotical stability 
with respect to initial values, 
and the dissipativity of solutions to nonlinear F-ODEs, we first recall some relevant results for the classical ODEs, namely,   
\begin{equation} \label{eq2}
\begin{split}
\frac{d}{dt}x(t)=f(t, x),~~~x\in \mathbb{R}^{d},
\end{split}
\end{equation}
which are assumed to have a unique solution $x\in C[[t_{0},+\infty),\mathbb{R}^{d}]$ for any given initial value $x(t_{0})=x_{0}$. 

In order to extend the concept of $A$-stability for linear multistep methods from the linear test equation to nonlinear systems, 
Dahlquist \cite{Dahlquist75} introduced the one-sided Lipschitz condition in 1975 for the ODEs (\ref{eq2}):
\begin{equation} \label{eq3}
\begin{split}
\langle f(t, x)-f(t, y), ~  x-y\rangle \leq \lambda  \|x-y\|^2,~\hbox{for all}~x,y\in  \mathbb{R}^{d},
\end{split}
\end{equation}
where $\lambda$ is the one-sided Lipschitz constant,  $\langle \cdot,\cdot\rangle$ and $\|\cdot,\cdot\|$ are 
the standard inner product and norm in $\mathbb{R}^{d}$. 
Then any two solutions $x(t)$ and $y(t)$ of equations (\ref{eq2}) with
different initial values $x_{0}$ and $y_{0}$ meet the following stability estimate:
\begin{equation} \label{eq4}
\begin{split}
\|x(t)-y(t)\|\leq\|x_{0}-y_{0}\|\cdot e^{\lambda(t-t_{0})}.
\end{split}
\end{equation} 
This implies the contractivity and exponential stability of the solutions to the ODEs (\ref{eq2}) 
with respect to the initial values for $\lambda\leq 0$ and $\lambda< 0$ respectively.

The one-sided Lipschitz condition (\ref{eq3}) has a significant influence on the numerical methods for 
stiff ODEs \cite{Butcher2, Hairer}. 
Stiff problems usually have large classical Lipschitz constant, but there may be a moderately sized, or even a negative one-sided Lipschitz constant. 
One class of important examples of stiff ODEs are derived from the space discretization of some parabolic equations such as reaction diffusion equations.
Dahlquist \cite{Dahlquist75} proposed the concept of $G$-stability for one-leg methods
and the corresponding linear multistep methods (LMMs) for stiff ODEs satisfying the one-sided Lipschitz condition. 
The fundamental equivalence between the $G$-stability and $A$-stability of LMMs and one-leg methods was 
established in 1978 \cite{Dahlquist}.  Moreover, Butcher \cite{Butcher} studied the
contractivity for Runge-Kutta methods and introduced the concept
of the $B$-stability; see the monograph \cite{Hairer} for more details.

Another type of ODE systems that are very close to the contractive ODEs 
is the so-called dissipative systems.  
The main feature of the dissipative systems is the presence of certain mechanisms of energy dissipation, 
which can lead to quite complicated limit regimes and structures  \cite{Hale}.  
For the ODEs  (\ref{eq2}), Humphries and Stuart \cite{Hump} imposed a structural condition on $f$, namely, 
\begin{equation} \label{eq5}
\begin{split}
\langle f(t, x), ~x\rangle \leq a-b\|x\|^2 \quad \hbox{for all}~x\in  \mathbb{R}^{d}
\end{split}
\end{equation}
for some $a\geq 0$ and $b>0$,  
which leads to the decay estimate of the form 
\begin{equation} \label{eq6}
\begin{split}
\|x(t)\|^2\leq \|x_{0}\|^2
e^{-2b(t-t_{0})}+\frac{a}{b}\left(1-e^{-2b(t-t_{0})}\right)\,.
\end{split}
\end{equation}
Hence the open ball $B(0,\sqrt{a/b}+\varepsilon)$ is an absorbing set as $t\to +\infty$ for
any given $\varepsilon>0$ and any given initial data.  
As defined in \cite{Hump}, an ODEs system is said to be dissipative if for any initial value $x_{0}$, 
there exists a time $t^*{(x_{0})}\ge t_{0}$ such that $x(t)\in B$ for $t>t^{*}$ and 
a bounded absorbing set $B$.
\footnote{
As noted in \cite{Hump}, 
a numerical method preserves the contractivity in (\ref{eq4}) is sometimes referred to be dissipative in the numerical literature,
but this conflicts with the corresponding terminology in dynamical systems.  
At the same time, the authors in \cite{Hump} give an accurate definition of dissipativity for ODEs, 
which mainly emphasizes the existence of a global attracting set. 
This definition was later widely accepted in the numerical literature, and we also follow this definition.}
We can easily see the exponential stability of $x(t)$ directly from (\ref{eq6}) for $a=0, b>0$.
It is known that a dissipative system should satisfy the one-sided Lipschitz condition.

There are various models of dissipative differential equations from physics and engineering;
see \cite{Hale,Temam}.
In 1994, Humphries and Stuart \cite{Hump} first studied the
numerical dissipativity for Runge-Kutta methods. They proved that
for DJ-irreducible Runge-Kutta methods, the algebraic stability is
sufficient to imply the dissipativity of the numerical solutions to 
(\ref{eq2}) with the dissipative condition (\ref{eq5}). Based on Dahlquist's $G$-stability theory
\cite{Dahlquist}, Hill \cite{Hill} demonstrated that the $A$-stability is equivalent to the dissipativity of LMMs and one-leg
methods for ODEs with the condition (\ref{eq5}).

It is very interesting and natural for us to understand if the fundamental results we have reviewed above 
about contractivity and dissipativity 
of the classical ODEs \eqref{eq2} can be established also for F-ODEs.
We first studied the Caputo F-ODEs in \cite{WangD} 
and established the contractivity and dissipativity under the same conditions as those for classical ODEs. 
More precisely, we obtained the following results  \cite{WangD}.

\begin{lemma} \label{lem:lem1}
\emph{(i)} Under the one-sided Lipschitz condition (\ref{eq3}) on $f$, it holds 
for any two solutions $x(t)$ and $y(t)$ to the F-ODEs \eqref{eq1} 
with two initial values $x_{0}$ and $y_{0}$ that 
\begin{equation}\label{eq7}
\begin{split}
\|x(t)-y(t)\|^{2}\leq \|x_{0}-y_{0}\|^{2}\cdot E_{\alpha}(2\lambda t^{\alpha}).
\end{split}
\end{equation}
In particular, we have that $\|x(t)-y(t)\|\leq
\|x_{0}-y_{0}\|$ for $\lambda\leq 0$.

\emph{(ii)} Let $x(t)$ be the solution of the F-ODEs (\ref{eq1}) and 
$f$ satisfy the dissipative condition (\ref{eq5}), then the fractional order system is dissipative 
in the sense that 
\begin{equation}\label{eq8}
\begin{split}
\|x(t)\|^2\leq \|x_{0}\|^2 E_{\alpha}\left[(-2b)
t^{\alpha}\right]+2a\int_{0}^{t}\frac{1}{(t-\tau)^{1-\alpha}}E_{\alpha,\alpha}[(-2b)(t-\tau)^{\alpha}]d\tau\,.
\end{split}
\end{equation}
Clearly, for any given $\varepsilon>0$, the ball $B(0,\sqrt{a/b}+\varepsilon)$ is
an absorbing set as $t\to +\infty$. 
\end{lemma}

To continue our discussions, we introduce two important functions, namely 
the Mittag-Leffler function $E_{\alpha}(z)$
and the generalized Mittag-Leffler function $E_{\alpha,\beta}(z)$ defined for $z\in \mathbb{C}$:
\begin{equation*}
\begin{split}
E_{\alpha}(z)=\sum_{k=0}^{\infty}\frac{z^{k}}{\Gamma(\alpha
k+1)},~\alpha>0;~~~E_{\alpha,\beta}(z)=\sum_{k=0}^{\infty}\frac{z^{k}}{\Gamma(\alpha
k+\beta)},~\alpha,\beta>0,
\end{split}
\end{equation*}
which are the fractional generalization of the exponential function
and play an important role in fractional calculus.
For $\alpha\in (0,1)$, these two functions have the following nice properties \cite{Kil, Pod}:
\begin{equation}\label{eq9}
\begin{split}
E_{\alpha}(t)=E_{\alpha,1}(t)>0,~~E_{\alpha,\alpha}(t)>0,~~\frac{d}{d
t}E_{\alpha,\alpha}(t)>0.
\end{split}
\end{equation}
By means of the asymptotic expansion of the Mittag-Leffler function \cite{Kil, Pod}, 
\begin{equation*}
\begin{split}
E_{\alpha,\beta}(\lambda t)=-\sum_{k=1}^{N}\frac{1}{\Gamma(\beta-k\alpha)}\frac{1}{(\lambda t)^{k}}+O\left(\frac{1}{(\lambda t)^{N+1}}\right)
\end{split}
\end{equation*}
for $N\in\mathbb{N}^{+}$, $t\rightarrow +\infty$ and $\lambda<0$,
we can obtain an explicit contractivity and dissipativity rates from (\ref{eq7}) and (\ref{eq8}), namely, 
it holds for some $c_{\alpha}>0$, 
\begin{equation}\label{eq10}
\begin{split}
\|x(t)-y(t)\|^{2}\leq \|x_{0}-y_{0}\|^{2}\cdot\frac{c_{\alpha}}{t^{\alpha}},~ \hbox{as}~ t \rightarrow +\infty\,, 
\end{split}
\end{equation}
\begin{equation}\label{eq11}
\begin{split}
\|x(t)\|^{2}\leq \|x_{0}\|^{2}\cdot\frac{c_{\alpha}}{t^{\alpha}}+\frac{a}{b}, ~ \hbox{as}~ t \rightarrow +\infty\,. 
\end{split}
\end{equation}
Here and in the rest of this work, we use $c_{\alpha}$ to represent a generic positive constant, 
which may take different values at different occasions, depending on $\alpha$ 
but independent of time $t$ or discrete time points $n$.

We may readily observe from \eqref{eq4},  (\ref{eq6}), \eqref{eq10} and (\ref{eq11}) 
that the contractivity and dissipativity rates with regard to initial values are exponential for ODEs 
while they are algebraic for F-ODEs. This reflects 
an essential difference between the long-term decay rates of solutions to classical initial value problems 
and fractional ones, mainly due to 
the nonlocal nature of fractional derivatives in some sense.

In another recent work \cite{WangZou}, we further studied the long-time stability of the solutions to stiff nonlinear fractional  functional differential equations (F-FDEs)
by means of a novel fractional delay-dependent Halanary-type inequality.
We investigated in \cite{WangZou} the effects of various functional terms such as time delay and delay integro-differential terms on the long-term properties of solutions.
A variety of complex dynamic behaviors were observed for the solutions to 
F-FDEs due to the involvement of functional terms and fractional derivatives.
In particular, we demonstrated rigorously 
the accurate algebraic decay rate the solutions observe 
with respect to various complex function perturbations in a given initial range.

In view of structure-preserving numerical methods, it
is desirable that the numerical solutions can inherit the long time behavior of the solutions to time fractional evolution equations.
This motivates one of the main focuses of this paper, namely, 
to study the contractivity and dissipativity of solutions to the numerical F-BDFs for nonlinear F-ODEs.
As we shall demonstrate both analytically and numerically, it is quite remarkable 
that the numerical solutions preserve exactly the same algebraic contractivity and dissipativity rates as the ones 
their continuous counterparts possess, described in (\ref{eq10}) and (\ref{eq11}).

We like to emphasize that contractivity and dissipativity for time fractional evolution equations are 
stronger decay behaviors than the usual stability.
Contractivity and dissipativity preserving numerical methods are more effective and desired 
in applications
than those stable schemes without such long-time characteristics,
especially when the solutions have various discontinuous points. 
Based on the two important lemmas established in this paper, 
we constructed in \cite{WangZou}  two effective difference schemes for F-FDEs, 
and proved that their numerical solutions preserve exactly the same algebraic contractivity rate as the 
one the continuous solutions observe.
The key idea in establishing this algebraic contractivity rate 
was to control various functional items through several new techniques.

To the best of our knowledge, 
the existing numerical stability analysis of F-ODEs is mostly focused on linear problems \cite{CermakJ, Garrappa2, Garrappa1,Garrappa3, Lubich1}, 
or nonlinear problems based on the classical Lipschitz hypothesis \cite{CaoW1, JinB3}.  In particular, Cao et al. proposed the time splitting schemes in \cite{CaoW1} 
and implicit-explicit difference schemes in \cite{CaoW2} to deal with stiff nonlinear F-ODEs.  
The methods in  \cite{CaoW1, CaoW2} have good linear stability without nonlinear iterations, 
but the special structures of the nonlinear function $f$ were not discussed. 
Noting that the nonlinear F-ODEs (\ref{eq1}) can be written equivalently 
as the Abel-Volterra integral equations of second kind with weakly singular kernel, 
the long time behavior was studied in \cite{Egger, Nevan} for 
the numerical solutions of the corresponding integral equations.
The error estimates were also obtained in \cite{Egger, Nevan}, under some stronger conditions
on the function $f$, requiring simultaneously the monotone condition (close to one-sided Lipschitz condition) and the global Lipschitz condition.

The rest of the paper is organized as follows.  
In \cref{sec:ImpT}, the contractivity rate obtained in \cite{WangD} is improved for scalar F-ODEs based on nonnegative preserving properties of the solution.
This result is then used to establish optimal numerical contractivity rate of F-BDFs for scalar F-ODEs in \cref{sec:Impmain}.
In \cref{sec:BDFs},  the contractivity and dissipativity of the numerical solutions to F-BDFs are established.
In \cref{sec:Prel}, a discrete version of the fractional generalization of the Leibniz rule is first obtained, 
which allows us to derive an energy-type inequality.   
Then a new asymptotical behavior is studied for the solution to a linear Volterra difference equation 
with algebraic decay rate,
which leads to the long time algebraic decay rate of the solutions to F-BDFs. The main results is proved in  \cref{sec:main}, 
and two typical examples of F-BDFs based on Gr\"{u}nwald-Letnikov and L1 difference schemes are presented in \cref{sec:exam}. 
The contractivity and dissipativity of some high order numerical schemes are 
developed in \cref{sec:HighO}, under slightly stronger conditions. 
Several numerical examples and the concluding remarks  are provided in \cref{sec:test} and \cref{sec:dis}, respectively.

\section{Improved contractivity rate of solutions to scalar F-ODEs}
\label{sec:ImpT}
In this section, we first derive a new contractivity rate of the solutions to the scalar F-ODE:
\begin{equation}\label{eq:scalar}
^{C}_{0}D_{t}^{\alpha}x(t)=f(x), \quad t>0,  ~ x\in \mathbb{R}\,, 
\end{equation}
under the one-sided Lipschitz condition (\ref{eq5}).
This improves the main results in \cite{WangD} and 
can be applied directly to establish optimal contractivity rate of numerical solutions 
to \eqref{eq:scalar} in \cref{sec:Impmain}, and 
to the spatial semi-discrete model of linear fractional sub-diffusion equation in \cref{sec:test}. 
The main tool in the analysis is 
the nonnegative preserving properties of the solutions to F-ODEs under appropriate conditions.

If the  F-ODEs (\ref{eq1}) is linear and stable, i.e., $f(x)=Ax$, where $A$ is a constant coefficient matrix, 
then we know the contractivity rate $\|x(t)-y(t)\|=O(t^{-\alpha})$ from the basic stability theory \cite{Matignon}.
But the rate was shown to become slower for general nonlinear  F-ODEs \cite{WangD}, 
namely, $\|x(t)-y(t)\|=O(t^{-\alpha/2})$; see  (\ref{eq10}). 
The energy analysis was used in \cite{WangD} to estimate the decay rate of $\|x(t)-y(t)\|^2$, 
which is bounded by $E_{\alpha}(2\lambda t^\alpha)$. However, we do not have 
$\sqrt{E_{\alpha}(2\lambda t^\alpha)}= E_{\alpha}(\lambda t^\alpha)$ for the Mittag-Leffler function, unlike 
the identity $\sqrt{e^{2\lambda t}}=e^{\lambda t}$  for the classical exponential function.
This is the main reason that causes the slower decay rate by the analysis in \cite{WangD}. 

We now make use of a new analytical tool to improve the above result to the optimal contractivity rate, 
namely, $\|x(t)-y(t)\|=O(t^{-\alpha})$ for nonlinear scalar F-ODE. 
The basic idea is to estimate the decay rate of $\|x(t)-y(t)\|$ directly, not $\|x(t)-y(t)\|^2$ as it did in \cite{WangD}.
This enables us to avoid the square-root operation of the Mittag-Leffler function.
To do this, we first present two auxiliary results.

\begin{lemma}[\cite{JinB2}]
\label{lem:JinBT}
For any
$x\in C[0, T ]\cap C^{1}(0, T ]$, if $x(t)$ attains its minimum at 
$t_{1}\in (0,T]$, 
then $^{C}_{0}D_{t_{1}}^{\alpha}x(t_{1})\leq 0$.
\end{lemma}

\begin{lemma}\label{lem:Pos}
(i) Under the dissipation condition \eqref{eq5} with $a=0$, 
if $x$ is a solution to the equation \eqref{eq:scalar}
and $x\in C[0, +\infty)\cap C^{1}(0, +\infty)$, then a positive initial value 
$x(0)$ implies $x(t)\ge 0$ for all $t>0$.

(ii) Under the one-sided Lipschitz condition (\ref{eq3}) on $f$ for some $\lambda <0$,
if $x$ and $y$ are two solutions to the equation \eqref{eq:scalar} such that 
$x, y\in C[0, +\infty)\cap C^{1}(0, +\infty)$
and $x(0)>y(0)$, then $x(t)\ge y(t)$ for all $t>0$.

\end{lemma}
\begin{proof}
We prove by contradiction. Assume there exists a time $t_{1}\in (0, T]$ for some $T>0$ such that $x(t_{1})<0$. 
Then we can find a time $t_{2}\in (0, T]$ such that $x(t_{2})=\min\limits_{t\in(0, T]} x(t)<0$, 
hence we know $~^{C}_{0}D_{t_{2}}^{\alpha}x(t_{2})\leq 0$ 
from \cref{lem:JinBT}. Using this result, we derive 
\begin{equation}\label{eq12}
\begin{split}
0\leq \left\langle ~^{C}_{0}D_{t_{2}}^{\alpha}x(t_{2}), x(t_{2})\right\rangle = \langle f(x(t_{2})), x(t_{2})\rangle \leq \lambda \|x(t_{2})\|^2<0.
\end{split}
\end{equation}
This contradiction yields the desired result in (i). The result in (ii) can be proved by the same argument.
\end{proof}

\begin{theorem}\label{thm:ImpT}
Under the same conditions on $f$, $x$ and $y$ as in Lemma \ref{lem:Pos}(ii) except that 
$x(0)$ may not be bigger than $y(0)$, the following asymptotic estimate holds
\begin{equation}\label{eq13}
\begin{split}
\|x(t)-y(t)\|\leq \|x_{0}-y_{0}\|\cdot\frac{c_{\alpha}}{t^{\alpha}} ~~ \hbox{as}~ t \rightarrow +\infty. 
\end{split}
\end{equation}\end{theorem}

\begin{proof}
Let $z(t)=x(t)-y(t)$, and we assume $z(0)>0$ (the same argument for $z(0)<0$). 
We readily see $z(t)\ge 0$ for all $t>0$ from \cref{lem:Pos} by noting that 
\begin{equation}\label{eq14}
\begin{split}
 \left\langle ~^{C}_{0}D_{t}^{\alpha} z(t), z(t)\right\rangle = \left\langle f(x(t))-f(y(t)), z(t)\right\rangle \leq \lambda \|z(t)\|^2<0\,.
\end{split}
\end{equation}
It follows also from (\ref{eq14}) that  $~^{C}_{0}D_{t}^{\alpha} z(t)\leq \lambda \|z(t)\|= \lambda z(t)$, which yields that 
$z(t)\leq z(0) E_{\alpha}(\lambda t^{\alpha})$. 
Now the desired estimate follows from the asymptotic expansion of the Mittag-Leffler function. 
\end{proof}

By a similar argument to the one of Theorem \ref{thm:ImpT} above, we can show that 
the dissipativity rate in (\ref{eq11}) for the solution $x(t)$ to the scalar F-ODE  \eqref{eq:scalar} 
can be improved: 
\begin{equation}\label{eq15}
\begin{split}
\|x(t)\|\leq \|x_{0}\|\cdot\frac{c_{\alpha}}{t^{\alpha}} ~~ \hbox{as}~ t \rightarrow +\infty, 
\end{split}
\end{equation}
under the dissipativity condition (\ref{eq5}) with $a=0$ and $b>0$.

For general vector-valued functions $z\in\mathbb{R}^{d}$ with $d>1$, 
we cannot expect to derive  similar results above from the inequality (\ref{eq14}). 
But we guess that the contractivity rate obtained in (\ref{eq10}) should be also 
optimal for nonlinear systems that can not be decoupled 
by diagonalization.  We will study this again in the numerical experiments. 

\section{Contractivity and dissipativity analysis of F-BDFs}
\label{sec:BDFs}
In this section we investigate the contractivity and dissipativity of numerical solutions to F-BDFs. 
For this, we introduce a step-size parameter $h>0$ and the corresponding time nodal points 
$t_{n}=nh, n=0,1,2,3,\cdots$. Further, we write $x_{n}$ for the approximation of  $x(t_{n})$ 
and $f_{n}=f(t_{n}, x_{n})$. 
As in integer-order differential equations, one basic approach of constructing difference schemes is based on the
numerical differentiation of fractional derivatives.
Because of the nonlocal nature of fractional derivatives, the
numerical approximation involves all discrete time points from $t_{0}$ to
$t_{n}$, leading to the numerical method for F-ODEs (\ref{eq1}) in the following full-term recursion
\begin{equation}\label{eq16}
\begin{split}
\sum_{j=0}^{n}\omega_{n-j}x_{j}=h^{\alpha}f(t_{n}, x_{n}),~~~n=1,2,3,\cdots.
\end{split}
\end{equation}
There are several approaches in the literature for determining the weight coefficients $\{\omega_{n}\}_{n=0}^{\infty}$, 
which yield a wide variety of numerical methods with different accuracies and stabilities. 
Solutions to time fractional equations often exhibit weak singularities at the origin, 
resulting in slower convergence rates of numerical solutions. 
Correction formulas were developed in \cite{Lubich2, JinLZ} to restore the convergence rate.
Since the correction terms do not affect the stability of numerical methods, 
we shall not consider them in this work.

Due to the major characteristic difference between F-ODEs and ODEs, 
the traditional analytical tools developed by Dahlquist \cite{Dahlquist} can not easily extended to F-ODEs.
It is well known that the concept of $G$-stability plays a
central role \cite{Dahlquist, Hairer, Hill} in the study of the contractivity and dissipativity of LMMs and one-leg
methods for classical ODEs, and 
most analyses are performed under the $G$-norm in $\mathbb{R}^{d\cdot k}$.
Unfortunately, the $G$-norm can not extend to F-LMMs, 
mainly because of the nonlocal nature of fractional operators, for which 
the dimension of the $G$ matrix increases with the time and is no longer fixed.

\subsection{Preliminaries}\label{sec:Prel}
In this subsection, we present some auxiliary results for the subsequent analysis. 
An inner product inequality involving Caputo fractional derivatives played a key
role in our analysis of F-ODEs \cite{WangD}, and the inequality was originated from 
the following important equality by Alikhanov, which is a fractional variant of the classical Leibniz formula.
\begin{lemma}[\cite{Alik}]\label{lem:Alikhanov}
For any two absolutely continuous functions $x(t)$ and $y(t)$ on $[0,T]$, 
the following equality holds for $0<\alpha<1$:
\beqn
&& \qq x^{T}(t)\cdot~^{C}_{0}D_{t}^{\alpha}y(t) +y^{T}(t)\cdot~^{C}_{0}D_{t}^{\alpha}x(t)= \label{eq:leibniz}\\
&&~^{C}_{0}D_{t}^{\alpha}\left(x^{T}(t)\cdot y(t)\right)+\frac{\alpha}{\Gamma(1-\alpha)}
\int_{0}^{t}\frac{1}{(t-\xi)^{1-\alpha}}\left(\int_{0}^{\xi}\frac{x'(\eta)d\eta}{(t-\eta)^{\alpha}}\cdot
\int_{0}^{\xi}\frac{y'(s)ds}{(t-s)^{\alpha}}\right)d\xi. \nb
\eqn 
\end{lemma}
We can easily derive the inequality 
\begin{equation}\label{eq18}
~^{C}_{0}D_{t}^{\alpha}\left(x^{T}(t)\cdot x(t)\right)\leq 2x^{T}(t)\cdot~^{C}_{0}D_{t}^{\alpha}x(t)~~for~~0<\alpha<1
\end{equation}
by taking $x(t)=y(t)$ in the identity \eqref{eq:leibniz} and noting the fact that 
$
\frac{\alpha}{\Gamma(1-\alpha)}
\int_{0}^{t}\frac{d\xi}{(t-\xi)^{1-\alpha}}\left(\int_{0}^{\xi}\frac{x'(\eta)d\eta}{(t-\eta)^{\alpha}}\right)^{2}\geq 0.
$

We shall start with the contractivity and dissipativity of numerical solutions to F-BDF (\ref{eq16}) 
under the following general assumptions on the weights $\{\omega_{n}\}_{n=0}^{\infty}$:
\begin{equation*}
\begin{split}
Assumption ~(A):~~~\left\{
    \begin{array}{l}
      \hbox{(i)} ~\omega_{0}>0, \\
      \hbox{(ii)}~ \omega_{j}\leq 0~\hbox{for~all}~j\geq 1, \\
      \hbox{(iii)}~ \sum\limits_{j=0}^n\omega_{j}\geq 0~\hbox{for~any~given}~n\geq 1\,, \\
    \end{array}
  \right.
\end{split}
\end{equation*}
then apply the results to two specific F-BDFs, based on Gr\"{u}nwald-Letnikov formula and L1 method.

As it is seen, the main motivation of Assumption $(A)$ is for deriving the following 
discrete version of the inequality (\ref{eq18}),
which is crucial to help us establish the numerical dissipativity and contractivity of F-BDFs.

\begin{lemma}\label{lem:Alikhanov2}
Under Assumption $(A)$, it holds for the F-BDF (\ref{eq16}):
\begin{equation}\label{eq19}
\begin{split}
\sum_{j=0}^{n}\omega_{n-j}\|x_{j}\|^2 \leq \left\langle 2x_{n}, \sum_{j=0}^{n}\omega_{n-j}x_{j}\right\rangle,~~~~n\geq 1.
\end{split}
\end{equation}
\end{lemma}

\begin{proof} 
The desired result comes from the direct calculations:
\begin{equation*}
\begin{split}
&\left\langle 2x_{n}, \sum_{j=0}^{n}\omega_{n-j}x_{j}\right\rangle-\sum_{j=0}^{n}\omega_{n-j}\|x_{j}\|^2\\
=&\left\langle 2x_{n}, \sum_{j=0}^{n}\omega_{n-j}x_{j}\right\rangle-\sum_{j=0}^{n}\omega_{n-j}\|x_{n}\|^2-\sum_{j=0}^{n}\omega_{n-j}\|x_{j}\|^2+\sum_{j=0}^{n}\omega_{j}\|x_{n}\|^2\\
=&\left\langle 2x_{n}, \sum_{j=0}^{n-1}\omega_{n-j}x_{j}\right\rangle-\sum_{j=0}^{n-1}\omega_{n-j}\|x_{n}\|^2-\sum_{j=0}^{n-1}\omega_{n-j}\|x_{j}\|^2+\sum_{j=0}^{n}\omega_{j}\|x_{n}\|^2\\
\geq&-\sum_{j=0}^{n-1}\omega_{n-j}\left(\left\|x_{n}\right\|-\left\|x_{j}\right\|\right)^2
+\left(\sum_{j=0}^{n}\omega_{j}\right)\|x_{n}\|^2\geq 0.
\end{split}
\end{equation*}
\end{proof}

The numerical discretization of F-ODEs often leads to some Volterra
difference equations of convolution type. The relevant results and
analytical tools for Volterra difference equations are often
employed to study the stability and asymptotic behaviors of
fractional numerical schemes \cite{CermakJ}. 
We now introduce some important results 
on the boundedness and asymptotic decay rate for the solutions to a class of linear convolution Volterra
difference equations. 

General speaking, it is much more difficult to achieve the exact decay rates of the solutions to difference equations 
than to establish qualitative properties such as the stability or asymptotic stability of some equilibrium solutions.
Applelby, Gy\H{o}ri and Rennolds  \cite{Applelby} derived exact convergence rates of some linear Volterra difference equations 
by making use of an elegant three-term decomposition of the discrete convolution.   
The relevant concepts and main results in \cite{Applelby} are included in \cref{sec:appendixA}.
A  remarkable advantage of \cref{lem:app} in \cref{sec:appendixA} is that the class of kernels could decay sub-exponentially, 
which allows us to derive the non-exponential convergence rates  for some asymptotically stable nontrivial solutions. 
This approach applies also to the difference schemes of F-ODEs, so 
we shall adopt it to derive the boundedness and exact contractivity rate of the F-BDF (\ref{eq16}). 

\begin{lemma}\label{lem:Rate}
Consider the Volterra difference equation 
\begin{equation}\label{eq20}
x_{n+1}= f_{n}+\sum\limits_{j=0}^{n}F_{n-j}x_{j},~~~n\geq1
\end{equation}
where the coefficients satisfy 
 $f_{n}\to \frac{c_{1}}{n^{\alpha}}, ~ F_{n}\to \frac{c_{2}}{n^{1+\alpha}}~as~ n\to \infty,~ and ~ \rho=\sum\limits_{j=0}^{\infty}|F_{j}|<1,$
for some constants $c_{1},  c_{2}>0$ and $0<\alpha<1$. 
Then we have the asymptotic estimate 
\begin{equation}\label{eq21}
 x_{n}\to  \frac{c_{1}\left(1-\rho \right)^{-1}}{n^{\alpha}}~ as~ n\to \infty.
\end{equation}
\end{lemma}

\begin{proof} 
We can not apply \cref{lem:app} in \cref{sec:appendixA} directly for the desired result
due to the facts that $f_{n}\to c_{1}/n^{\alpha}$ and $\gamma_{n}=1/(n+1)^{\alpha} \notin W(1)$ for $0<\alpha<1$, and that the series 
$\sum_{j=1}^{\infty}\frac{1}{n^{\alpha}}$ diverges.

We introduce a simple transformation $y_{n}=\frac{x_{n}}{n}$ for $n\ge 1$,
 and let $y_{0}=x_{0}$, $g_{n}=\frac{f_{n}}{n+1}$ and $G_{n, j}=\frac{j}{n+1} F_{n-j}$.  
Then the equation (\ref{eq20}) becomes 
\begin{equation} \label{eqa3}
\begin{split}
y_{n+1}=g_{n}+\sum_{j=0}^{n}G_{n, j}~ y_{i}, ~~n\ge 1.
\end{split}
\end{equation}

We may note that (\ref{eq20}) is
a convolution difference equation while  equation (\ref{eqa3}) is not. 
Obviously,  it holds that 
$g_{n}\to c_{1}/n^{1+\alpha}$ as  $n\to \infty$. Following the idea developed in \cite{Applelby}, we now take the weight sequence 
 $\gamma_{n}=1/(n+1)^{1+\alpha}\in W(1)$ and compute $L_{\gamma}(y)
 =\lim_{n\to \infty} {y_{n}}/{\gamma_{n}}$
 to give  a non-trivial limit, 
 which yields that $y_{n}$ behaves like $O(n^{-(1+\alpha)})$ asymptotically. 
Letting $z_{n}=y_{n} / \gamma_{n}$, we can rewrite equation  (\ref{eqa3}) as
\begin{equation} \label{eqa4}
\begin{split}
z_{n+1}=h_{n}+\sum_{j=0}^{n}H_{n, i}~ z_{j}, ~~n\ge 1,
\end{split}
\end{equation}
with 
$h_{n}=\frac{g_{n}}{\gamma _{n+1}},~~ H_{n, j}=\frac{\gamma_{j}}{\gamma_{n+1}}G_{n, j}=\frac{j}{n+1} \frac{\gamma_{j}}{\gamma_{n+1}}F_{n-j}.$
Now we plan to derive the limit of $z_{n}$ satisfying equation  (\ref{eqa4}) by \cref{lem:app}. 
It suffices to verify all the conditions in the lemma. 
Firstly, it is easy to see 
\begin{equation*} 
\begin{split}
\lim_{n\to \infty}\sup \sum_{j=0}^{m}|H_{n, n-j}|=\sum_{j=0}^{m}|F_{j}| \lim_{n\to \infty} \left(\frac{n-j}{n+1}\frac{\gamma_{n-j}}{\gamma_{n+1}}\right)=\sum_{j=0}^{m}|F_{j}|,
\end{split}
\end{equation*}
therefore, $\lim\limits_{m\to \infty}\sup \lim\limits_{n\to \infty}\sup \sum_{j=0}^{m}|H_{n, n-j}|=\sum_{j=0}^{\infty}| F_{j}|=\rho<1.$
Secondly, for any fixed $m>0$,  we have
\begin{equation*} 
\begin{split}
\lim_{n\to \infty}\sup \sum_{j=m}^{n-m}|H_{n, n-j}|
\le \sup_{j\ge 0}\frac{|F_{j}|}{\gamma_{j}}  \lim_{n\to \infty} \frac{\gamma_{n}}{\gamma_{n+1}} 
\lim_{n\to \infty} \left(\frac{n-j}{n+1}\frac{\gamma_{n-j}\gamma_{j}}{\gamma_{n+1}}\right),
\end{split}
\end{equation*}
which implies that $\lim\limits_{m\to \infty}\sup\left( \lim\limits_{n\to \infty}\sup \sum_{j=m}^{n-m}|H_{n, n-j}|\right)=0.$
It remains to compute that
\begin{equation*} 
\begin{split}
 \lim_{n\to \infty}H_{n, m}= \lim_{n\to \infty} \left(\frac{m}{n+1} \frac{\gamma_{n-m}}{\gamma_{n+1}} \frac{F_{n-m}}{\gamma_{n-m}} \right)\gamma_{m},
\end{split}
\end{equation*}
which implies that $H_{\infty, m}=\lim\limits_{n\to \infty}H_{n, m}=0$.
In fact, it follows directly from \cref{lem:app} that $\lim\limits_{n\to \infty}z_{n}= L_{\gamma}(y)= \left(1-\rho \right)^{-1} L_{\gamma}(g)$,
which leads readily to our desired estimate 
\begin{equation*}
 y_{n}=\frac{x_{n}}{n}  \to  \frac{c_{1}\left(1-\rho \right)^{-1}}{n^{1+\alpha}}~ \hbox{as}~ n\to \infty.
\end{equation*}
\end{proof}

The discrete energy inequality in \cref{lem:Alikhanov2} and the $O(n^{-\alpha})$ decay rate of Volterra difference equation in \cref{lem:Rate} are crucial in our subsequent analysis.
They are also very useful for analyzing the long-term stability and decay rate of other more complex problems, such as F-FDEs \cite{WangZou} and time fractional PDEs.

Although we consider only the uniform grids in this work, 
we can trace our whole analysis to find out that our results are also true for non-uniform grids 
as long as the corresponding weight coefficients meet the specified assumptions, 
including the popular graded grids and the non-uniform L1 formula \cite{LiaoH}.
This is very useful when we construct adaptive numerical methods or schemes with relatively large time steps.

\subsection{Numerical contractivity and dissipativity}\label{sec:main}
We now present one of our main results in this paper, which can be seen as  the discrete version of \cref{lem:lem1}.

\begin{theorem}\label{thm:thm1}
Assume that the weights $\{\omega_{k}\}_{n=0}^{\infty}$ of the F-BDF (\ref{eq16}) satisfy Assumption $(A)$,
and there exists a constant $c_{\alpha}>0$ such
that $|\omega_{k}|\leq c_{\alpha}/k^{1+\alpha}$ for $1\leq k\leq n-1$ and
$|\omega_{n}|\leq c_{\alpha}/n^{\alpha}$ for any $n\in N^{+}$.

\emph{(i)} If function $f$ in (\ref{eq16}) satisfies the one-sided Lipschitz condition \eqref{eq3}, 
and
$\rho_{1}=\sum_{j=1}^{\infty} \frac{|\omega_{j}|}{\omega_{0}-2\lambda h^{\alpha}}<1$ for any $h>0$,
then the F-BDF (\ref{eq16})  is contractive, and its solution can preserve the exact contractivity rate  
as the true solution to F-ODEs \eqref{eq1} (cf.\,(\ref{eq10})), namely, 
$\|x_{n}-y_{n}\|^2\leq c_{1} {\left\|x_{0}-y_{0}\right\|^2}{n^{-\alpha}}$ as $n\to \infty$, 
with $c_{1}=\left(1-\rho_{1} \right)^{-1}\frac{c_{\alpha}}{\omega_{0}-2\lambda h^{\alpha}}$.

\emph{(ii)} If function $f$ in (\ref{eq16}) satisfies condition \eqref{eq5},
and $\rho_{2}=\sum_{j=1}^{\infty}\frac{|\omega_{j}|}{\omega_{0}+2b h^{\alpha}}<1$ for any $h>0$, 
then the F-BDF is dissipative,  i.e., for any given initial value $x_{0}$ and $\varepsilon>0$, there is a bounded set
$B\left(0,r\right)$ and $n_{0}\in N^{+}$ such that $x_{n}\in B\left(0,r\right)$ for all $n\geq n_{0}$, 
with 
$r=\sqrt{{c_2a}/{b}}+\varepsilon$ and $c_{2}=\left(1- \rho_{2}\right)^{-1}$.
Moreover, if the condition \eqref{eq5} is satisfied with $a=0$, 
the numerical solution has the exact dissipativity rate as the exact solution 
to F-ODEs \eqref{eq1}, namely, 
$\|x_{n}\|^2\leq c_{3} {\left\|x_{0}\right\|^2}{n^{-\alpha}}$ as $n\to \infty$, 
with $c_{3}=\left(1-\rho_{2} \right)^{-1}\frac{c_{\alpha}}{\omega_{0}+2b h^{\alpha}}$. 

\end{theorem}

\begin{proof}
(i) Let $x_{j}$ and $y_{j}$ be the numerical solutions of the F-BDF (\ref{eq16})
with two different initial values $x_{0}$ and $y_{0}$, respectively. Put $z_{n}=x_{n}-y_{n}, n\geq 0$.
We can easily see that $\sum\limits_{j=0}^{n}\omega_{n-j}z_{j}=h^{\alpha}\left(f( x_{n})-f(y_{n})\right).$
Taking the inner product with $2z_{n}$ on both sides 
and applying the one-sided Lipschitz condition and  \cref{lem:Alikhanov2}, we get
\begin{equation}\label{eq22}
\begin{split}
\sum_{j=0}^{n}\omega_{n-j}\left\|z_{j}\right\|^{2}\leq 2\lambda h^{\alpha}
\left\|z_{n}\right\|^{2},
\end{split}
\end{equation}
which can be rewritten as $\left(\omega_{0}-2\lambda h^{\alpha}\right)
\left\|z_{n}\right\|^{2}\leq\sum\limits_{j=0}^{n-1}|\omega_{n-j}|\left\|z_{j}\right\|^{2}$ by noting that $\omega_{n-j}<0$ for $j=0,1,...,n-1$.
Since the weights $\omega_{n}$ and $\omega_{j}$ for $j\le n-1$ have different decay rates, the above Volterra difference inequality can be further rewritten 
\begin{equation*}
\begin{split}
\|z_{n}\|^2\leq \frac{|\omega_{n}|}{\omega_{0}-2\lambda h^{\alpha}}\left\|z_{0}\right\|^2+\sum\limits_{j=1}^{n-1}\frac{|\omega_{n-j}|}{\omega_{0}-2\lambda h^{\alpha}}\left\|z_{j}\right\|^2.
\end{split}
\end{equation*}
Now applying \cref{lem:Rate} yields the desired decay rate 
$\|z_{n}\|^2\leq \left(1-\rho_{1} \right)^{-1}\frac{\left\|z_{0}\right\|^2}{\omega_{0}-2\lambda h^{\alpha}} \frac{c_{\alpha}}{n^{\alpha}}=c_{1}  \frac{\left\|z_{0}\right\|^2}{n^{\alpha}}$ as $n\to \infty$.\\

(ii) It follows directly from the dissipativity condition (\ref{eq5}) and \cref{lem:Alikhanov2} that
\begin{equation*}
\begin{split}
\left\langle 2x_{n}, \sum_{j=0}^{n}\omega_{n-j}x_{j}\right\rangle=2h^{\alpha}\left\langle f_{n}, x_{n}\right\rangle
\leq 2h^{\alpha}\left(a-b\|x_{n}\|^2\right),\\
\end{split}
\end{equation*}
which implies that $\left(\omega_{0}+2h^{\alpha}b\right)\|x_{n}\|^2\leq 2h^{\alpha}a-\sum\limits_{j=0}^{n-1}\omega_{n-j}\left\|x_{j}\right\|^2$ for $n\geq1$, 
leading to the convolution Volterra inequality 
\begin{equation*}
\begin{split}
\|x_{n}\|^2 &\leq \frac{2h^{\alpha}a}{\omega_{0}+2h^{\alpha}b}+\sum\limits_{j=0}^{n-1}\frac{|\omega_{n-j}|}{\omega_{0}+2h^{\alpha}b}\left\|x_{j}\right\|^2~ for~n\geq1.
\end{split}
\end{equation*}
By applying \cref{lem:app} in \cref{sec:appendixA}, we obtain that $\|x_{n}\|^2 \leq \left(1-\rho_{2} \right)^{-1} \frac{2h^{\alpha}a}{\omega_{0}+2h^{\alpha}b} \leq c_{2}\frac{a}{b}$ as $n\to \infty$.
 where $\rho_{2}=\sum_{j=1}^{\infty}\frac{|\omega_{j}|}{\omega_{0}+2 h^{\alpha}b}<1$ and $c_{2}=\left(1-\rho_{2} \right)^{-1}$.
 The poof of the desired dissipativity rate for $a=0, b>0$ is similar to the proof of (i) and omitted here.
\end{proof}

One may observe from the proof of part (ii) in \cref{thm:thm1},  it is not easy to derive 
the dissipativity result by the usual Gr\"{o}nwall-like inequalities, 
because those estimates depend often directly on the initial values $x_{0}$, 
 but the dissipativity is a long time feature of solutions to F-ODEs and is independent of the initial values. 
 
An alternative approach for the boundedness of $\|x_{n}\|$ in  part (ii) of  \cref{thm:thm1} 
is to apply some discrete variants of a Paley-Wiener theorem, which was introduced by Lubich  \cite{Lubich1}.
We demonstrate below that the results obtained by this approach is consistent to the ones 
in \cref{thm:thm1}.
We first recall a result from \cite{Lubich1}.

\begin{lemma} \label{lem:Lubich}
Consider the discrete Volterra  equation $y_{n}= p_{n}+\sum_{j=0}^{n}q_{n-j}y_{j}, ~n\geq0$,
 where the kernel $\{q_{n}\}_{n=0}^{\infty}$ belongs to $l^{1}$, i.e.,   $\sum_{j=0}^{\infty}|q_{j}|<\infty$.
Then $y_{n}\to 0$  (resp. bounded) whenever $p_{n}\to 0$ (resp. bounded) as $n\to\infty$
if and only if the Paley-Wiener condition is satisfied, i.e.,
\begin{equation}\label{eq24}
\sum\limits_{j=0}^{\infty}q_{j}\zeta^{j}\ne 1 ~for~ |\zeta|\le 1.
\end{equation}
\end{lemma}

If we define a sequence $\{r_{n}\}_{n=0}^{\infty}$ by $\frac{1}{1-\sum\limits_{j=0}^{\infty}q_{j}\zeta^{j}}=\sum\limits_{j=0}^{\infty}r_{j}\zeta^{j}$,
we can easily check from the proof of \cref{lem:Lubich} that 
if $\{q_{n}\}_{n=0}^{\infty}$ belongs to $l^{1}$ and the Paley-Wiener condition (\ref{eq24}) holds, 
then $\{r_{n}\}_{n=0}^{\infty}$ is  also in $l^{1}$, 
and the estimate holds $\|y\|_{l^{\infty}}\le\|r\|_{l^{1}}\|p\|_{l^{\infty}}.$

Now consider the Volterra difference equation related to  (ii) of  \cref{thm:thm1}, i.e.,
\begin{equation*}
\begin{split}
\|x_{n}\|^2 = \frac{2h^{\alpha}a}{\omega_{0}+2h^{\alpha}b}+\sum\limits_{j=0}^{n-1}\frac{|\omega_{n-j}|}{\omega_{0}+2h^{\alpha}b}\left\|x_{j}\right\|^2~ for~n\ge 1.
\end{split}
\end{equation*}
The assumption $\rho_{2}=\sum\limits_{j=1}^{\infty}\frac{|\omega_{j}|}{\omega_{0}+2b h^{\alpha}}<1$ implies that 
the kernel $\left\{\frac{|\omega_{j}|}{\omega_{0}+2b h^{\alpha}}\right\}_{j=1}^{\infty}$ belongs to $l^{1}$ and
the corresponding Paley-Wiener condition 
$\sum\limits_{j=1}^{\infty}\frac{|\omega_{j}|\zeta^{j}}{\omega_{0}+2b h^{\alpha}}\ne 1$ for $|\zeta|\le 1$
holds. Let 
\begin{equation}\label{eq27}
\frac{1}{1- \sum\limits_{j=1}^{\infty}\frac{|\omega_{j}|\zeta^{j}}{\omega_{0}+2b h^{\alpha}}}=\sum\limits_{j=0}^{\infty}r_{j}\zeta^{j} ~for~ |\zeta|\le 1.
\end{equation}
Then $ r_{j}\ge 0$ and $\{r_{n}\}_{n=0}^{\infty}$ is in $l^{1}$.
Taking $\zeta=1$ in (\ref{eq27}) yields that 
 \begin{equation*}
 \|x_{n}\|^2\le \|r\|_{l^{1}} \frac{2h^{\alpha}a}{\omega_{0}+2h^{\alpha}b}=   \left(1- \rho_{2}\right) ^{-1} \frac{2h^{\alpha}a}{\omega_{0}+2h^{\alpha}b}
 \le \left(1- \rho_{2}\right) ^{-1} \frac{a}{b} ~~ \hbox{as}~ n\to\infty,
\end{equation*}
which is the same as the corressponding results in part (ii) of \cref{thm:thm1}.

\subsection{Improved numerical contractivity rates for scalar F-ODEs}\label{sec:Impmain}
In \cref{sec:ImpT}, we presented an optimal contractivity rate 
for the scalar F-ODE (\ref{eq:scalar}).
A typical application of this new result is for the spatial semi-discrete model of linear fractional sub-diffusion equation, 
and this will be carefully validated by numerical experiments in \cref{sec:test}. 
Next we demonstrate that this optimal contractivity rate can be preserved exactly by the numerical solutions to 
the scalar F-ODE (\ref{eq:scalar}).
We first derive some nonnegative preserving properties of the numerical solutions to the F-BDF 
(\ref{eq16}) for (\ref{eq:scalar}). 

\begin{lemma}\label{lem:NumPP}
Let function $f$ in the scalar F-ODE (\ref{eq:scalar}) satisfy that 
$\langle f(x), x\rangle \leq \lambda \|x\|^2$ for some $\lambda <0$, and 
$x(t)$ be the solution to (\ref{eq:scalar}) with $x(0)>0$. 
Then under Assumption (A), the solutions to the F-BDF (\ref{eq16}) are all nonnegative. 
\end{lemma}

\begin{proof}
We prove by mathematical induction.
For $n=1$, we see directly from (\ref{eq16}) that 
$\omega_{1}x_{0}+\omega_{0}x_{1}=h^{\alpha}f(x_{1})$. Then we can get 
by taking the inner product with $x_{1}$ on both sides that $-\omega_{1}x_{0}x_{1}\geq(\omega_{0}-h^{\alpha}\lambda)x_{1}^2\geq 0$, 
which implies $x_{1}\geq 0$.

We now prove $x_{n}\geq 0$ under the condition that $x_{j}\geq 0$ for $j=1,2,...,n-1$. 
Taking the inner product with $x_{n}$ in both sides of the F-BDF (\ref{eq16}) gives that 
$\langle x_{n}, \sum_{j=0}^{n}\omega_{n-j}x_{j}\rangle=h^{\alpha}\left\langle f_{n}, x_{n}\right\rangle\leq h^{\alpha}\lambda\|x_{n}\|^2,$
which can be rewritten as 
\begin{equation*}
\begin{split}
 x_{n}\left(-\sum_{j=0}^{n-1}\omega_{n-j}x_{j}\right)\geq (\omega_{0}-h^{\alpha}\lambda)\|x_{n}\|^2\geq 0.\\
\end{split}
\end{equation*}
This implies that $x_{n}\geq 0$.
\end{proof}

\begin{theorem}\label{thm:ImpN}
Let function $f$ in the scalar F-ODE (\ref{eq:scalar}) satisfy the one-sided Lipschitz condition (\ref{eq3}), 
$x_n$ and $y_n$ are two solutions to the F-BDF (\ref{eq16}) with different initial values $x_0$ and $y_0$. 
Then under Assumption (A), the following contractivity estimate holds
\begin{equation}\label{eq28}
\begin{split}
\|x_{n}-y_{n}\|\leq \|x_{0}-y_{0}\|\cdot\frac{c_{1}}{n^{\alpha}},~ \hbox{as}~ n \rightarrow \infty, 
\end{split}
\end{equation}
where $c_{1}$ is the same as in \cref{thm:thm1}.

On the other hand, 
if function $f$ in the scalar F-ODE (\ref{eq:scalar}) 
satisfies the dissipative condition (\ref{eq5}) with $a=0$, 
then the solutions to the F-BDF (\ref{eq16}) decay as 
\begin{equation}\label{eq29}
\begin{split}
\|x_{n}\|\leq \|x_{0}\|\cdot\frac{c_{3}}{n^{\alpha}} ~ ~\hbox{as}~ n \rightarrow \infty \q \hbox{($c_{3}$ is the same 
as in \cref{thm:thm1})}. 
\end{split}
\end{equation}
\end{theorem}

\begin{proof}
The proof is very similar to the one of \cref{thm:thm1},
but we estimate the decay rate of $z_{n}=x_{n}-y_{n}$ directly rather than $\|x_{n}-y_{n}\|^2$,
which allows us to avoid the square-root operation of the Mittag-Leffler function.
Without lose of generality, we assume $z_{0}=x_{0}-y_{0}>0$. Using the dissipative condition (\ref{eq5}), 
we can derive 
\begin{equation}\label{eq30}
\begin{split}
\Big\langle z_{n}, \sum\limits_{j=0}^{n}\omega_{n-j}z_{j}\Big\rangle=h^{\alpha}\left\langle  z_{n}, f(x_{n})-f(y_{n})\right\rangle\leq h^{\alpha}\lambda\|z_{n}\|^2,
\end{split}
\end{equation}
then it follows from \cref{lem:NumPP} that $z_{n}\ge 0$ for all $n\geq1$.  
This non-negativity and  the inequality  (\ref{eq30}) yield that
$\sum_{j=0}^{n}\omega_{n-j}z_{j}\leq h^{\alpha}\lambda \|z_{n}\|= h^{\alpha}\lambda z_{n}$, from 
which and \cref{lem:Rate} the contractivity rate (\ref{eq28}) follows readily. 

The proof of the dissipativity rate (\ref{eq29}) can be done similarly.
\end{proof}

\subsection{Examples of F-BDFs}\label{sec:exam}
In this subsecion, we present some concrete examples of F-BDFs, whose 
weights $\{\omega_{j}\}_{j=0}^{\infty}$ meet all the conditions required for the results 
we have derived in the previous three subsections.
We consider two widely used low-order schemes of the form F-BDF (\ref{eq16}), 
i.e., the Gr\"{u}nwald-Letnikov formula \cite{Pod} and the  L1 method \cite{LinY, SunZ}.
The coefficients of these two schemes meet very nice properties so  
we can establish the energy-type inequality in \cref{lem:Alikhanov2} and the decay rate in  \cref{lem:Rate} 
directly. 

\subsubsection{Gr\"{u}nwald-Letnikov formula}
\label{sec:GL}

The wildly used Gr\"{u}nwald-Letnikov (G-L) fractional derivative \cite{Kil,Pod} are defined for $0<\alpha<1$ by
\begin{equation}\label{eq32}
\begin{split}
~^{GL}_{~0}{D}^\alpha_tx(t)=\lim_{h\rightarrow 0^{+}}\frac{\left(\Delta_{h}^{\alpha}\right)x(t)}{h^{\alpha}}=
\lim_{ h\rightarrow 0^{+}}\frac{1}{h^{\alpha}}\sum_{k=0}^{m=[t/h]}(-1)^{k}\left(
                                                                                 \begin{array}{c}
                                                                                   \alpha \\
                                                                                   k \\
                                                                                 \end{array}
                                                                               \right)
                                                                              x(t-kh).
\end{split}
\end{equation}

If we do not perform the limit operation $h\rightarrow 0^{+}$ in \eqref{eq32} but take $h>0$ to be the step-size,
then the discretized version of the operator $~^{GL}_{~0}{D}^\alpha_tx(t)$ can be expressed as
\begin{equation}\label{eq:gl}
\begin{split}
~^{GL}_{~0}{D}^\alpha_tx(t_{n})=
\frac{1}{h^{\alpha}}\sum_{k=0}^{n}(-1)^{k}\left(
                                                                                 \begin{array}{c}
                                                                                   \alpha \\
                                                                                   k \\
                                                                                 \end{array}
                                                                               \right)
                                                                              x_{n-k}+O(h)
                                                                              =\frac{1}{h^{\alpha}}\sum_{k=0}^{n}\omega_{k}x_{n-k}+O(h),
\end{split}
\end{equation}
where the coefficients are given by $\omega_{k}=(-1)^{k}\left(
                                                                                 \begin{array}{c}
                                                                                   \alpha \\
                                                                                   k \\
                                                                                 \end{array}
                                                                               \right)$, $k=0, 1, \cdots, n$. 
For the Caputo derivative, we introduce the following scheme \cite{Garrappa2}:
 \begin{equation}\label{eq34}
\begin{split}
~^{C}_0{D}^\alpha_{t_{n}}x(t_{n})=\frac{1}{h^{\alpha}} \Big(\sum\limits_{j=1}^{n}\omega_{n-j}x_{j}+\delta_{n}x_{0}\Big)+O(h),
\end{split}
\end{equation}
where the coefficient $\delta_{n}$ is set to be $\delta_{n}=-\sum_{j=0}^{n-1}\omega_{j}$ so that 
the sum of the weights in \eqref{eq34} equals to zero, 
which is beneficial to the numerical stability of the scheme \cite{Garrappa2}. 
The G-L formula is a simple and effective numerical scheme with first order accuracy, 
and its weights $\omega_{k}$ meet the following properties. 
\begin{lemma}[\cite{Garrappa2}]\label{lem:GL}
For $0<\alpha<1$, the coefficients $\omega_{k}=(-1)^{k}\left(
                                                                                 \begin{array}{c}
                                                                                   \alpha \\
                                                                                   k \\
                                                                                 \end{array}
                                                                               \right)$ 
satisfy                   
                                                         
\emph{(i)} $\omega_{0}=1$, $\omega_{n}<0, |\omega_{n+1}|<|\omega_{n}|, ~~n=1,2,\cdots$;

\emph{(ii)} $\omega_{0}=-\sum\limits_{j=1}^{\infty}\omega_{j}>-\sum\limits_{j=1}^{n}\omega_{j}, ~~n\ge 1$;

\emph{(iii)} $\omega_{n}=O(n^{-1-\alpha}), ~\delta_{n}=O(n^{-\alpha})~as~n\to \infty$.

\end{lemma}

\subsubsection{L1 method}
\label{sec:L1}
The L1 method is among the most popular algorithms for the discretization of the Caputo derivative.
It often leads to unconditionally stable algorithms, and has the accuracy $O(h^{2-\alpha})$ for smooth data \cite{LinY,SunZ} while has order $O(h)$ for non-smooth data in uniform grids \cite{LiaoH}.
The L1 method can be written as
\begin{equation}\label{eq:scheme}
{\small 
\begin{split}
^{C}_0{D}^\alpha_tx(t)|_{t=t_{n}}
=\frac{1}{h^{\alpha}}\sum_{k=0}^{n-1}b_{n-k-1}\left(x_{k+1}-x_{k}\right)+O(h^{q})
=\frac{1}{h^{\alpha}}\sum_{k=0}^{n}\gamma_{n-k}x_{k}+O(h^{q}),
\end{split}
}
\end{equation}
where the coefficients are given by 
$b_{k}=\frac{1}{\Gamma(2-\alpha)}\left((k+1)^{1-\alpha}-k^{1-\alpha}\right)$
for $0\le k\le n-1$, $\gamma_{0}=\frac{1}{\Gamma(2-\alpha)}$,
 $\gamma_{k}=\frac{1}{\Gamma(2-\alpha)}\left((k+1)^{1-\alpha}-2k^{1-\alpha}+(k-1)^{1-\alpha}\right)$  for  $k=1,2,...,n-1$, 
 and $\gamma_{n}=\frac{1}{\Gamma(2-\alpha)}\left((n-1)^{1-\alpha}-n^{1-\alpha}\right)$.
The second formula in \eqref{eq:scheme} is the discrete convolution quadrature and 
its coefficients have the following properties, which can be checked directly.
\begin{lemma}\label{lem:L1}
The coefficients of the L1 method meet the properties:
\begin{equation*}
\begin{split}
\emph{(i)} ~ & ~\gamma_{0}>0, \gamma_{1}<\gamma_{2}<\cdots<\gamma_{n-1}<0,~ \gamma_{n}<0 ~for~any~~ n\geq1; \\
\emph{(ii)} ~ &~ k^{1+\alpha}\gamma_{k}\rightarrow -\frac{\alpha}{\Gamma(1-\alpha)},~as~ k\rightarrow \infty~ for ~k\neq n,\\
~ & and~n^{\alpha}\gamma_{n}\rightarrow -\frac{1}{\Gamma(1-\alpha)},~as~ n\rightarrow \infty.
\end{split}
\end{equation*}
\end{lemma}

A common feature of the G-L formula and L1 method is that the sign of the weights $\{\omega_{j}\}_{j\geq 1}$ 
and $\{\gamma_{j}\}_{j\geq 1}$ remain negative, which are crucial to the results in \cref{lem:Alikhanov2}.  
But this feature is no long true for high order schemes, 
such as the fractional trapezoidal rule, the second order F-BDF formula and fractional Newton-Gregory formula, 
and there are always some positive weights \cite{Garrappa3}.  
Another important feature of the G-L formula and L1 method is that their weights 
$\delta_{n}$ and $\gamma_{n}$ decay in the order $O(n^{-\alpha})$.
But the coefficients of the schemes for 
the F-ODEs with the Riemann-Liouville fractional derivative decay faster, namely, 
in the order $O(n^{-1-\alpha})$.
We emphasize that the decay rates of the weights $\delta_{n}$ and $\gamma_{n}$ essentially determine 
the decay rates of the numerical method (\ref{eq16}); see \cref{lem:Rate}.

For both the G-L formula and L1 method,  we now verify the conditions in \cref{thm:thm1} are satisfied. 
Indeed, for the G-L formula, we have
\begin{equation*}
\begin{split}
\rho_{1}&=\sum\limits_{j=1}^{\infty}\frac{|\omega_{j}|}{\omega_{0}-2\lambda h^{\alpha}}=\frac{1}{1-2\lambda h^{\alpha}}<1,
 ~~c_{1}= \left(1-\frac{1}{2\lambda h^{\alpha}}\right)\frac{c_{\alpha}}{1-2\lambda h^{\alpha}},\\
\rho_{2}&=\sum\limits_{j=1}^{\infty}\frac{|\omega_{j}|}{\omega_{0}+2b h^{\alpha}}=\frac{1}{1+2b h^{\alpha}}<1,~~~ c_{2}= 1+\frac{1}{2b h^{\alpha}}, 
\end{split}
\end{equation*}
while for  the L1 method, we have
\begin{equation*}
\begin{split}
\rho_{1}&=\frac{1}{1-2\Gamma(2-\alpha)\lambda h^{\alpha}}<1, 
~c_{1}= \left(1-\frac{1}{2\Gamma(2-\alpha)\lambda h^{\alpha}}\right)\frac{c_{\alpha}}{1-2\Gamma(2-\alpha)\lambda h^{\alpha}},\\
\rho_{2}&=\frac{1}{1+2\Gamma(2-\alpha)b h^{\alpha}}<1,~c_{2}= 1+\frac{1}{2b\Gamma(2-\alpha) h^{\alpha}}.
\end{split}
\end{equation*}
The above shows both methods satisfy all the conditions of \cref{thm:thm1}, therefore 
are contractive, dissipative, and preserve the optimal contractivity rate.

\subsection{High order numerical approximations}\label{sec:HighO}
In this subsection, we establish the contractivity and dissipativity of the F-BDFs 
for some high order approximations under slightly stronger conditions than those in (\ref{eq3}) and (\ref{eq5}). 
As two typical examples, we consider the second order F-BDFs \cite{Lubich2, Garrappa3} 
and a $3-\alpha$ order approximation based on quadratic interpolation approximations \cite{GaoG, LvX}.

\subsubsection{Second order F-BDFs}
The numerical method in (\ref{eq16}) can be written 
 \begin{equation*}
\begin{split}
~^{C}_0{D}^\alpha_{t_{n}}x(t_{n})=\frac{1}{h^{\alpha}} \Big(\sum_{j=1}^{n}\mu_{n-j}x_{j}+\delta_{n}x_{0}\Big)+O(h^{q}),
\end{split}
\end{equation*}
where the coefficient $\delta_{n}$ is given by $\delta_{n}=-\sum_{j=0}^{n-1}\mu_{j}$, 
while the weights $\{\mu_j\}$ are generated by the function 
 \begin{equation} \label{eq37}
\begin{split}
\mu(\xi)=\left(\frac{3}{2}-2\xi+\frac{1}{2}\xi^2\right)^{\alpha}
=\left(\frac{3}{2}\right)^{\alpha}(1-\xi)^{\alpha}\left(1-\frac{1}{3}\xi\right)^{\alpha}
=\sum\limits_{n=0}^{\infty}\mu_{j}\xi^{j}\,,
\end{split}
\end{equation}
and can be computed by $\mu_{j}=\left(\frac{3}{2}\right)^{\alpha}\sum_{l=0}^{j}3^{-l}\omega_{l}\omega_{j-l}$
\cite{ChenM},
where $\omega_{l}$ 
are the coefficients of the Gr\"{u}nwald-Letnikov formula
in \eqref{eq:gl}.
Using this formula and the asymptotic expansion of the binomial coefficients \cite{Lubich2}, 
we have the following more details about the behaviors of these coefficients, 
the asymptotic decay rate $\mu_{n}=O(n^{-\alpha-1})$ \cite{Lubich2} and 
$\delta_{n}=O(n^{-\alpha})$.
\begin{lemma}[\cite{ChenM}]\label{lem:BDF2} 
For $0<\alpha<1$,  we have
 \begin{equation*}
\begin{split}
&\mu_{0}=\left(\frac{3}{2}\right)^{\alpha}, ~\mu_{1}=-\left(\frac{3}{2}\right)^{\alpha}\frac{4\alpha}{3}, ~\mu_{2}=\left(\frac{3}{2}\right)^{\alpha}\frac{\alpha(8\alpha-5)}{9},\\
&\mu_{3}=\left(\frac{3}{2}\right)^{\alpha}\frac{4\alpha(\alpha-1)(7-8\alpha)}{81}; ~ \mu_{j}<0~ for~ j\geq 4;~~ \sum_{j=0}^{\infty}\mu_{j}=0; \\
&\mu_{n}=O(n^{-\alpha-1}) ~~\mbox{as} ~n\to \infty; \q \delta_{n}=O(n^{-\alpha})\,.
\end{split}
\end{equation*}
\end{lemma}

\subsubsection{Quadratic interpolation approximations}
The L1 method can be seen as a linear interpolation formula on each subinterval, and 
has the accuracy of order $2-\alpha$ for smooth functions. 
High order approximations to the Caputo derivative can be constructed by using the multiple nodal 
interpolations for the integrands. 
In particular, the quadratic interpolation approximation (QIA) \cite{GaoG, LvX} gives 
\begin{equation} \label{eq38}
\begin{split}
^{C}_0{D}^\alpha_tx(t)|_{t=t_{n}}=\frac{1}{h^{\alpha}}\sum_{j=0}^{n}\mu_{n-j}^{(n)}x_{j}+O(h^{3-\alpha}).
\end{split}
\end{equation} 
The coefficients $\{\mu_{j}^{(n)}\}_{j\geq 0}$ were given in  \cite{GaoG, LvX}, and have the following properties.

\begin{lemma}[\cite{LvX}] \label{lem:QIA}
For $n\geq 4$, $0<\alpha<1$ and $d_{0}=1/\Gamma(3-\alpha)$, we have
 \begin{equation*}
\begin{split}
&\mu_{0}^{(n)}=2^{1-\alpha}(1+\alpha/2)d_{0}>0, ~-\frac{4}{3}d_{0}<\mu_{1}^{(n)}<0, 
~-\frac{1}{3}d_{0}<\mu_{2}^{(n)}<\frac{1}{2}d_{0},\\
&\mu_{j}^{(n)}<0~ for~ j\geq 3; ~~\sum_{j=0}^{n}\mu_{j}^{(n)}=0;  ~~\mu_{n}^{(n)}=O(n^{-\alpha}) ~~as ~~n\to \infty.\\
\end{split}
\end{equation*}
\end{lemma}

As we can see from \cref{lem:BDF2} and \cref{lem:QIA}, 
some weights in the above two methods might be positive, 
e.g., $\mu_{2}$ and $\mu_{3}$ in second order F-BDFs, and $\mu_{2}^{(n)}$ in QIA.
This is very different from the previous  G-L formula and L1 method,
and causes some difficulties for us to derive the energy-like inequality in \cref{lem:Alikhanov2}.  
In \cite{LvX}, a special technique was introduced to transform all the coefficients $\{\mu_{j}^{(n)}\}_{j\geq 1}$ to be negative. 
But the transformation is essentially linear and appears to be difficult to apply to nonlinear systems.
However, we are still able to establish the contractivity and dissipativity for these high order schemes 
by slightly relaxing our previous assumptions, as it is shown in \cref{thm:HO}.

\begin{theorem}\label{thm:HO}
(i) Let function $f$ in the F-ODEs (\ref{eq1}) satisfy the one-sided Lipschitz condition (\ref{eq3}).
Then the second order F-BDFs and quadratic interpolation approximation are contractive if 
\begin{equation}\label{eq39}
\begin{split}
h^{\alpha}\lambda \leq
\left\{
\begin{aligned}
&-2\mu_{2}-2\mu_{3}~~\mbox{for second order F-BDFs},\\
&-2\mu_{2}^{(n)}~~\mbox{for the quadratic interpolation approximation},\\
\end{aligned}
\right.
\end{split}
\end{equation}
and any two different solutions $x_n$ and $y_n$ meet the following contractivity estimate 
\begin{equation}\label{eq40}
\begin{split}
\|x_{n}-y_{n}\|^2\leq \|x_{0}-y_{0}\|^2\cdot\frac{c_{\alpha}}{n^{\alpha}} ~~\hbox{as}~ n \rightarrow \infty. 
\end{split}
\end{equation}
(ii) 
Let function $f$ in (\ref{eq1}) satisfy the dissipative condition (\ref{eq5}).
Then the second order F-BDFs and quadratic interpolation approximation are dissipative if 
\begin{equation}\label{eq41}
\begin{split}
h^{\alpha}b \geq
\left\{
\begin{aligned}
&2\mu_{2}+2\mu_{3} ~~\mbox{for  second order F-BDFs},\\
&2\mu_{2}^{(n)}~~\mbox{for the quadratic interpolation approximation},\\
\end{aligned}
\right.
\end{split}
\end{equation}
 i.e., for any initial value $x_{0}$ and $\varepsilon>0$, there is a bounded set
$B\left(0,r\right)$ and $n_{0}\in N^{+}$ such that $x_{n}\in B\left(0,r\right)$ for all $n\geq n_{0}$, with 
$r=\sqrt{c_{\alpha} a/{b}}+\varepsilon$.
Moreover, if condition (\ref{eq5}) holds with $a=0$, the dissipativity is given by 
$\|x_{n}\|^2\leq c_{\alpha}{\left\|x_{0}\right\|^2}{n^{-\alpha}}$ as $n\to \infty$. 
\end{theorem}

\begin{proof}
(i) We prove only the result for the  QIA, and the same argument can be used to show the result for 
the second order F-BDFs.
Let $z_{n}=x_{n}-y_{n}$, and 
take the inner product with $2z_{n}$ on both sides of the numerical scheme 
and then apply the one-sided Lipschitz condition to obtain 
\begin{equation}\label{eq42}
\begin{split}
\left\langle \sum\limits_{j=0}^{n}\mu_{n-j}^{(n)}z_{j}, 2z_{n}\right\rangle=2h^{\alpha}\left\langle f(x_{n})-f(y_{n}), z_{n}\right\rangle\leq 2h^{\alpha}\lambda\|z_{n}\|^2.
\end{split}
\end{equation}

Without loss of generality, we assume $n\geq 4$.
We know from \cref{lem:QIA} that only the coefficient $\mu_{2}^{(n)}$ may be positive.
If $\mu_{2}^{(n)}\leq 0$, the results follow as in the proof of \cref{thm:thm1}. 
We now prove for the case that $\mu_{2}^{(n)}> 0$.
Define the new weights $\{\tilde{\mu}_{j}\}_{j\geq 0}$: 
$\tilde{\mu}_{0}=\mu_{0}^{(n)}+\mu_{2}^{(n)}, ~~ \tilde{\mu}_{1}=\mu_{1}^{(n)}, ~~\tilde{\mu}_{2}=0~ and ~\tilde{\mu}_{j}=\mu_{j}^{(n)}~ for ~j\geq 3.$
Then the inequality (\ref{eq42}) can be written as 
\begin{equation}\label{eq44}
\begin{split}
\left\langle \sum\limits_{j=0}^{n}\tilde{\mu}_{n-j}z_{j}, 2z_{n}\right\rangle+ \mu_{2}^{(n)}\left\langle z_{n-2}-z_{n}, 2z_{n}\right\rangle\leq 2h^{\alpha}\lambda\|z_{n}\|^2.
\end{split}
\end{equation} 

We can check that $\tilde{\mu}_{0}>0, ~ \tilde{\mu}_{j}\leq 0$ for all $j\geq 1$ and that $\sum\limits_{j=0}^{n}\tilde{\mu}_{j}\geq 0$ for $n\geq 1$.
Now \cref{lem:Alikhanov2} and the Cauchy inequality yields that 
$\sum\limits_{j=0}^{n}\tilde{\mu}_{n-j}\|z_{j}\|^{2} -3\mu_{2}^{(n)}\|z_{n}\|^{2}-\mu_{2}^{(n)}\|z_{n-2}\|^{2} \leq 2h^{\alpha}\lambda\|z_{n}\|^2.$
This inequality is equivalent to
\begin{equation}\label{eq46}
\begin{split}
\|z_{n}\|^{2}\leq \frac{1}{\tilde{\mu}_{0}-3\mu_{2}^{(n)}-2h^{\alpha}\lambda} \Big(\sum\limits_{j=0}^{n-1}| \tilde{\mu}_{n-j} | \|z_{j}\|^{2} +\mu_{2}^{(n)}\|z_{n-2}\|^{2}\Big).
\end{split}
\end{equation} 
Now the assumption $h^{\alpha}\lambda<-2\mu_{2}^{(n)}$ ensures that 
\begin{equation}\label{eq47}
\begin{split}
\rho_{3}= \frac{1}{\tilde{\mu}_{0}-3\mu_{2}^{(n)}-2h^{\alpha}\lambda} \Big(\sum\limits_{j=0}^{n-1}| \tilde{\mu}_{n-j} | +\mu_{2}^{(n)}\Big)
=\frac{\mu_{0}^{(n)}+2\mu_{2}^{(n)}}{\mu_{0}^{(n)}-2\mu_{2}^{(n)}-2h^{\alpha}\lambda} <1.
\end{split}
\end{equation} 
Then the estimate (\ref{eq40}) follows from \cref{lem:Rate}.

(ii) The same as in part (i), we introduce the new weights $\{\tilde{\mu}_{j}\}_{j\geq 0}$ and can then derive using 
the dissipative condition, 
\begin{equation}\label{eq48}
\begin{split}
\left\langle \sum\limits_{j=0}^{n}\tilde{\mu}_{n-j}x_{j}, 2x_{n}\right\rangle+ \mu_{2}^{(n)}\left\langle x_{n-2}-x_{n}, 2x_{n}\right\rangle\leq 2h^{\alpha}\left( a-b\|x_{n}\|^2\right).
\end{split}
\end{equation} 
By \cref{lem:Alikhanov2} and the Cauchy inequality, we can further deduce 
\begin{equation}\label{eq49}
\begin{split}
\|x_{n}\|^{2}\leq \frac{1}{\tilde{\mu}_{0}-3\mu_{2}^{(n)}+2h^{\alpha}b } \Big(2h^{\alpha}a+\sum\limits_{j=0}^{n-1}| \tilde{\mu}_{n-j} | \|x_{j}\|^{2} +\mu_{2}^{(n)}\|x_{n-2}\|^{2}\Big).
\end{split}
\end{equation} 
Now the assumption $h^{\alpha}b>2\mu_{2}^{(n)}$ ensures that 
\begin{equation}\label{eq50}
\begin{split}
\rho_{4}= \frac{1}{\tilde{\mu}_{0}-3\mu_{2}^{(n)}+2h^{\alpha}b} \Big(\sum\limits_{j=0}^{n-1}| \tilde{\mu}_{n-j} | +\mu_{2}^{(n)}\Big)
=\frac{\mu_{0}^{(n)}+2\mu_{2}^{(n)}}{\mu_{0}^{(n)}-2\mu_{2}^{(n)}+2h^{\alpha}b} <1.
\end{split}
\end{equation} 
Then we can see the desired dissipativity follows from \cref{lem:app} or the discrete Paley-Wiener theorem.
The dissipativity decay for $a=0$ can be derived similarly. 
\end{proof}

{\it Remark 3.1}. The conditions in (\ref{eq39}) and (\ref{eq41}) are automatically fulfilled when the weights of the numerical schemes are negative. 
This is the case with both the G-L formula and L1 method.
When the weight coefficients are positive, the constraints in (\ref{eq39}) and (\ref{eq41}) are not very restrictive, 
because the weights $\mu_{2}^{(n)}$ or $\mu_{2}+\mu_{3}$ are usually small, 
no more than $0.5$ which can be verified by simple calculation. 
So these conditions are relatively easy to meet.
We note that these constraints are only sufficient, and we still do not know if they are also necessary.

{\it Remark 3.2}.
For many other high order F-BDFs in the literature, such as Lubich's LMMs \cite{Lubich2}, forth order methods \cite{ChenM},  and higher order methods based on
interpolation formulas \cite{CaoJ}, etc., 
their weights $\{\mu_{j}\}_{j=0}^{\infty}$ often satisfy the conservation property $\sum_{j=0}^{\infty}\mu_{j}=0$ and $\mu_{0}>0$.
If there exist $k$ weights from $\{\mu_{j}\}_{j\geq 1}$ that are positive, say 
$\mu_{i_{1}}, \mu_{i_{2}}, ..., \mu_{i_{k}}$, we may naturally modify the condition in (\ref{eq39}) to 
$h^{\alpha}\lambda \leq -2(\mu_{i_{1}}+ \mu_{i_{2}}+...+\mu_{i_{k}}),$
then the contractivity of F-BDFs can be derived. 
The dissipativity can be also obtained in a similar manner.

\section{Numerical experiments}\label{sec:test}
In this section,  several numerical examples are presented to validate  
our theoretically predicted 
contractivity and dissipativity of the implicit F-BDFs (\ref{eq16}), 
and to reveal the algebraic decay rates of the F-ODEs. 
We will compare the numerical performance of F-BDFs with 
the popular predictor-corrector type methods, i.e.,  fractional Adams-Bashforth-Moulton (F-ABM) method, 
which was  proposed in \cite{Diethelm2}, especially for stiff problems. 

\subsection{Fractional Lorenz system} \label{sec:test1}
Consider the system
\begin{equation}\label{eqn1}
\begin{split}
\left\{
\begin{aligned}
~^{C}_{~0}D_{t}^{\alpha}x_{1}(t)&=x_{3}+(x_{2}-c_{1})x_{1},\\
~^{C}_{~0}D_{t}^{\alpha}x_{2}(t)&=1-c_{2}x_{2}-x_{1}^{2},\\
~^{C}_{~0}D_{t}^{\alpha}x_{3}(t)&=-x_{1}-c_{3}x_{3},\\
\end{aligned}
\right.
\end{split}
\end{equation}
where $c_{1}, c_{2}$ and $c_{3}$ are positive parameters and $c_{2}>1/2$.
This example contains many well-known dynamical systems such as the Lorenz, Chen, Chua systems and the financial system \cite{Pet}.
The classical Lorenz system was proved to be dissipative \cite{Hump} for $\alpha=1$.
Let $x=(x_{1},x_{2},x_{3})^{T}$, 
then we have by simple calculations that 
\begin{equation}\label{eqn2}
\begin{split}
\langle f(x),x\rangle&=-c_{1}x_{1}^{2}-c_{2}x_{2}^{2}-c_{3}x_{3}^{2}+x_{2}\\
&\leq\frac{1}{2}-c_{1}x_{1}^{2}-\left(c_{2}-\frac{1}{2}\right)x_{2}^{2}-c_{3}x_{3}^{2}\\
&\leq a-b\|x\|^{2}, 
\end{split}
\end{equation}
with $a=1/2$, $b=\min\{c_{1},c_{2}-1/2,c_{3}\}$. Thus the system is dissipative,
and the set $B(0,\sqrt{a/b}+\varepsilon)$ is absorbing. 

\begin{table}\label{tablea}
 \caption{Numerical performances of the F-ABM for Example 4.1} \label{tabEx1}
\begin{center}
\begin{tabular}{llllllllll }
\hline  & $\alpha=0.9$    & $\alpha=0.7$  & $\alpha=0.5$ & $\alpha=0.3$ & $\alpha=0.1$\\
\hline
Blowup  $h$ &5e-2   &2e-2         &4e-3         & 1e-4      & 2e-13  \\
Stable  $h$ &4e-2      &1e-2         &3e-3        & 5e-5       &1e-13    \\
\hline
\end{tabular}
\end{center}
\end{table}

\begin{figure}[H]
\begin{minipage}{1.0\linewidth}
\centering
\includegraphics[scale=0.4]{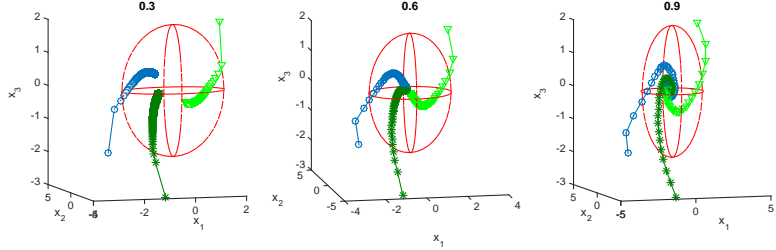}
\caption{Numerical solutions for $\alpha=0.3, 0.6$ and $0.9$ with parameters $c_{1}=1/4, c_{2}=1, c_{3}=1/4$, $a=1/2, b=1/4$ and the reference ball $B(0,\sqrt{2})$.  
Three orbits are computed by G-L method with $h=0.2$ and $T=100$ form initial values $(2, 1, 2)^{T},  (-2, 3, -2)^{T}$ and $(-1, -4, -3)^{T}$, respectively.} 
\label{ex11}
\end{minipage}
\end{figure}

\begin{figure}[H]
\begin{minipage}{1.0\linewidth}
\centering
\includegraphics[scale=0.45]{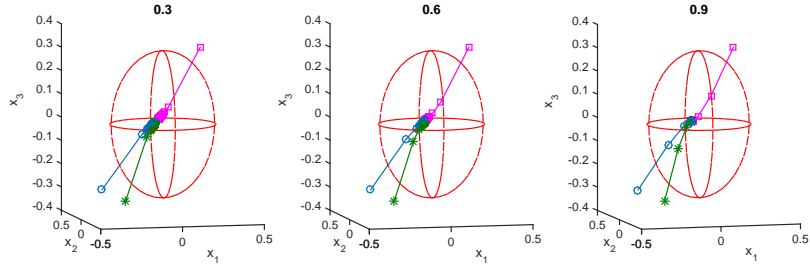}
\caption{Numerical solutions for $\alpha=0.3, 0.6$ and $0.9$ with parameters $c_{1}=5, c_{2}=6, c_{3}=5$, $a=1/2, b=5$ and the reference ball $B(0,1/\sqrt{10})$.
Three orbits are computed  by second order F-BDFs with $h=0.4$ and $T=200$ from the initial values $(0.3, 0.3, 0.3)^{T},  (-0.3, 0.3, -0.3)^{T}$ and $(-0.3, -0.3, -0.3)^{T}$, respectively.} 
\label{ex12}
\end{minipage}
\end{figure}

Figs.\ \ref{ex11} and \ref{ex12} plot the numerical solutions computed by the G-L formula  
and second F-BDFs respectively with various parameters and fractional order $\alpha$. 
They show that the order  $\alpha$ heavily affect the shape and size of the absorbing set, 
but all the computed solutions are kept inside the ball $B(0,\sqrt{a/b})$ when the time $t$ increases, as expected.
Comparing Fig.\,\ref{ex11} with Fig.\,\ref{ex12}, we can see that when $b$ is greater (i.e., the conditions in (\ref{eq41}) 
is satisfied), 
 the solution has stronger dissipation characteristics, which shrinks to the absorbing set $B(0, \sqrt{a/b})$ at a faster rate.
For the L1 method or QIA, similar numerical results are observed but not provided here.

In \cite{WangD}, the F-ABM method was employed to simulate this system.
In order to keep the stability,  the step size $h$ is required such that $h<h_{0}(\alpha)$ for some $h_{0}(\alpha)>0$.
Moreover, when the order $\alpha$ is small, this limitation usually becomes very demanding and can not be used for long time computation. 
As a comparison, we list in Tab.\ \ref{tablea} the step size limits that make the F-ABM method to be stable.
For $\alpha=0.1$, the step sizes have to be selected about $h=1e-13$, and the numerical blowup appears for $h=2e-13$.
But the F-BDFs method is stable uniformly for any $h>0$ and $\alpha \in (0, 1)$.
In fact, we guess that there exists certain equivalence relation between linear stability and numerical dissipativity for F-ODEs. 
Hill proved the corresponding equivalence theorem for classical ODEs in \cite{Hill}.

\subsection{Fractional sub-diffusion equation}\label{sec:test2}
Consider the 2D fractional sub-diffusion equation
\begin{equation}\label{eqn3}
\begin{split}
\left\{
\begin{aligned}
~^{C}_{~0}D_{t}^{\alpha}u(t,x,y)&=k\left(u_{xx}(t,x,y)+u_{yy}(t,x,y)\right)+g(t,x,y), ~(x,y)\in\Omega,\\
u(0,x,y)&=u_{0}(x,y),\\
u(t,x,y)&=0,~(x,y)\in\partial\Omega,
\end{aligned}
\right.
\end{split}
\end{equation}
where $\Omega=[0, 1]^2$ and the diffusion coefficient $k>0$.
Applying the standard finite element method with rectangular grids in the spatial direction, we get the F-ODEs 
\begin{equation}\label{eqn4}
\begin{split}
~^{C}_{~0}D_{t}^{\alpha}U(t)=-kAU(t)+G(t),
\end{split}
\end{equation}
where $U(t), G(t)\in R^{N_{x} \cdot N_{y}}$, and $N_{x}, N_{y}$
are the numbers of nodes in the $x, y$-directions respectively. 
It is well-known that the stiffness matrix $A$ is similar to a symmetric positive matrix, i.e., $D=P^{-1}AP$,  where $P$ is an orthogonal matrix.
Hence, its eigenvalues are positive and real, i.e., $0<\lambda_{1}\leq \lambda_{2} \leq\cdots\leq \lambda_{N_{x} \cdot N_{y}}$.
Let $F(U)=-kAU(t)+G(t)$. By direct calculations we have
\begin{equation}\label{eqn5}
\begin{split}
\langle F(U)-F(V),U-V\rangle=&-k(U-V)^{T}A(U-V)\\
=&-k(U-V)^{T}PDP^{-1}(U-V)\\
\leq & \mu\|U-V\|^2,
\end{split}
\end{equation}
where  $\mu=-k\lambda_{1}<0$. Therefore, the F-ODEs (\ref{eqn4}) satisfy the one-sided Lipschitz condition (\ref{eq3}), 
and they are contractive.
From (\ref{eq10}),  we have the contractivity rate  
\begin{equation}\label{eqn6}
\begin{split}
\|U(t)-V(t)\|^2\leq \|U(0)-V(0)\|^2\cdot \frac{c_{\alpha}}{t^{\alpha}}, ~~c_{\alpha}>0.
\end{split}
\end{equation}
 where $U(0), V(0)$ are two given initial values.
 Since $A$ is symmetric and positive in the semi-discrete system (\ref{eqn4}), it can be diagonalized, so the contractivity rate can be improved to be  
\begin{equation}\label{eqn6b}
\begin{split}
\|U(t)-V(t)\|\leq \|U(0)-V(0)\|\cdot \frac{c_{\alpha}}{t^{\alpha}}, ~~c_{\alpha}>0.
\end{split}
\end{equation}

\begin{figure}[!ht]
\begin{center}
$\begin{array}{cc}
 \includegraphics[scale=0.4]{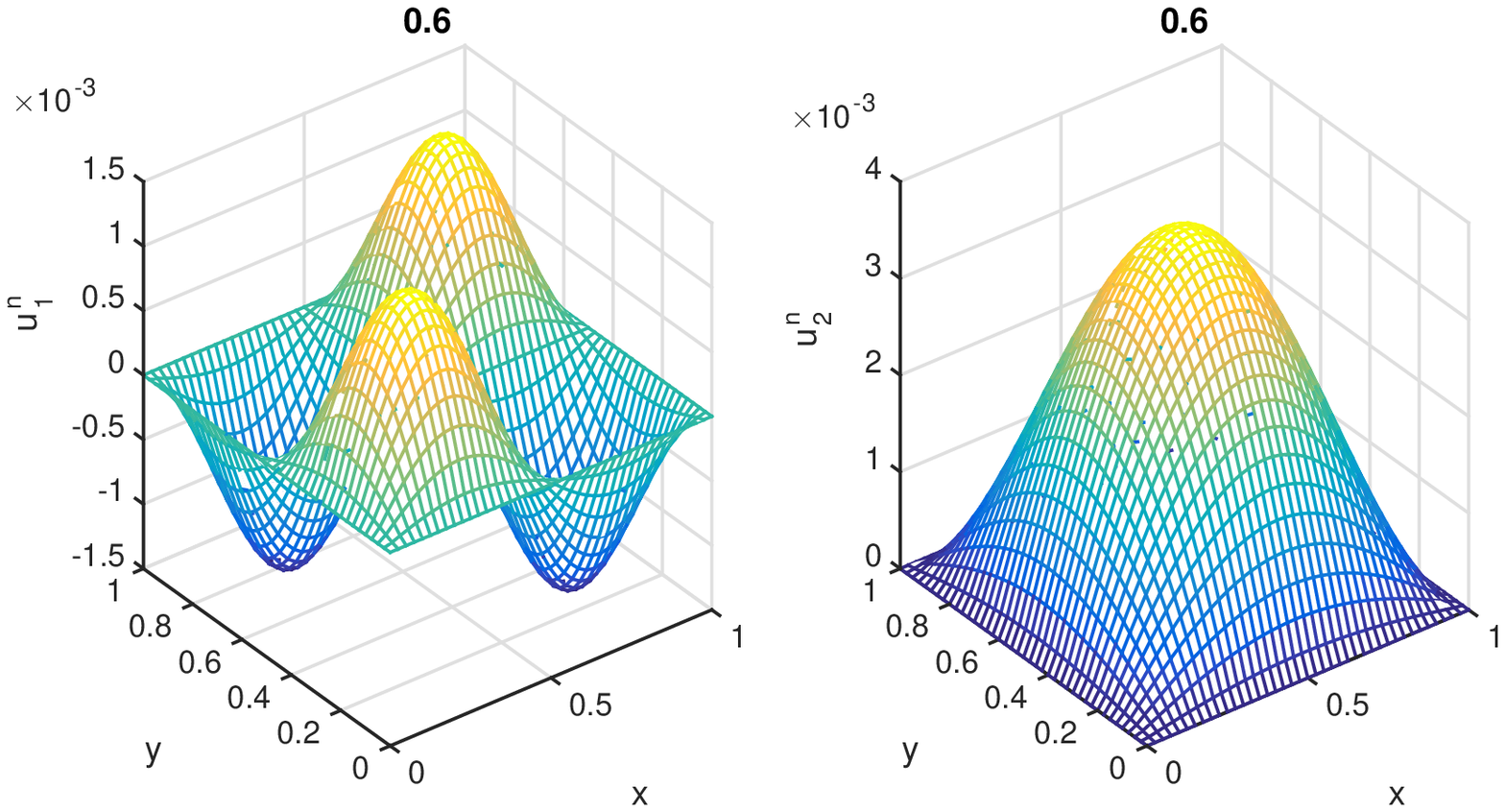}
 \includegraphics[scale=0.3]{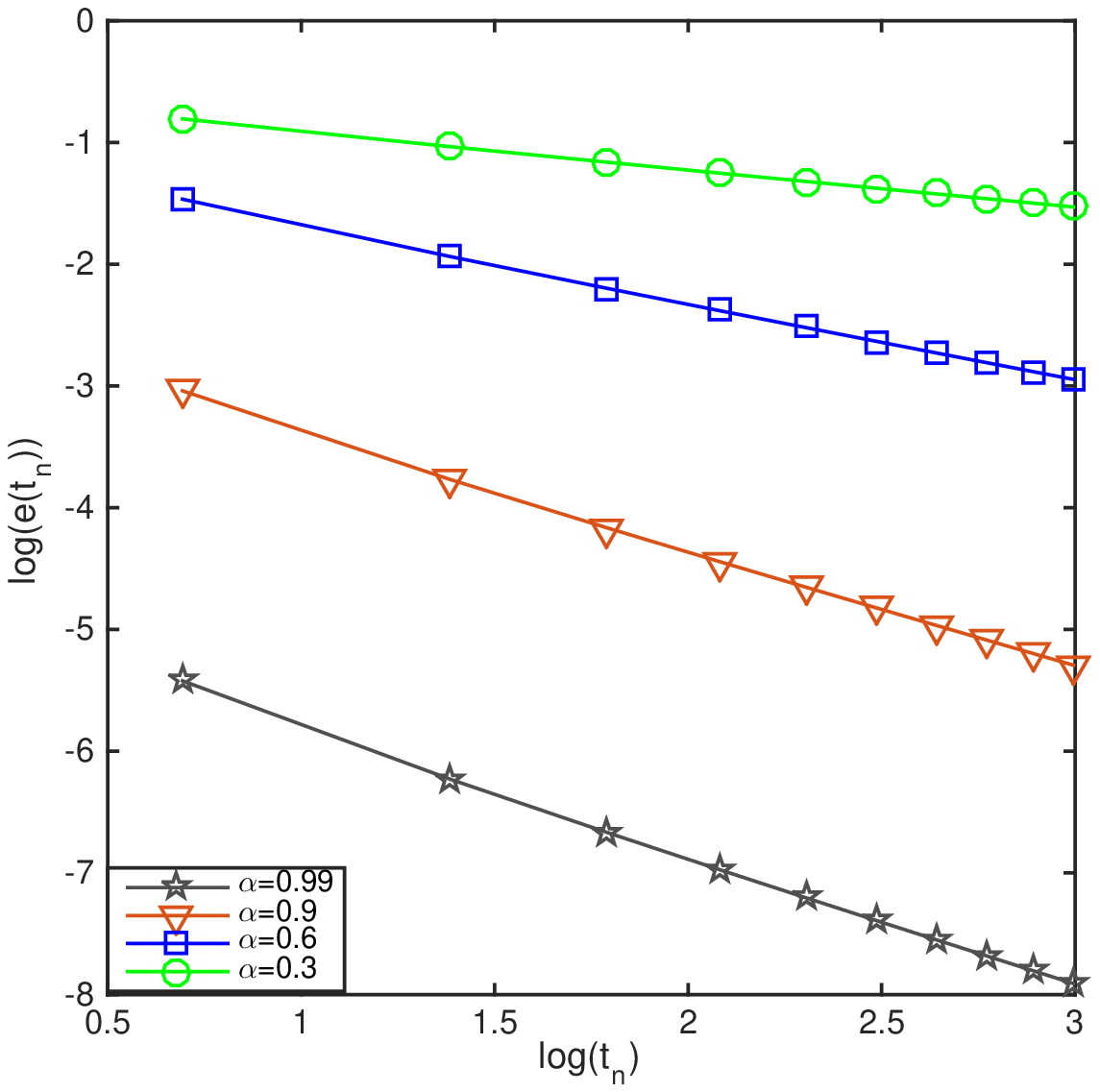}
 \end{array}$
\end{center}
\caption{Numerical solutions obtained by L1 method at $T=20$ for $h=0.2, \alpha=0.6$ with initial values $u_{1}$ and $u_{2}$, 
and the corresponding difference function $e(t)$ on [0, 100] for h=0.2 with $\alpha=0.3, 0.6, 0.9$ and $0.99$.}
\label{pic2}
\end{figure}

In the numerical simulation,  we take the initial values $u_{0}^{1}=\sin(2\pi x)\sin(2\pi y)$, $u_{0}^{2}=10xy(1-x)(1-y)$ and $g(t,x,y)=0$.
Let $e(t)=\|U(t)-V(t)\|$, then the discrete $l^2$-norm is given by 
$e(t_{n})=\Big(\frac{1}{N_{x} \cdot N_{y}}\sum_{k=1}^{N_{x} \cdot N_{y}} |U^{n}_{k}-V^{n}_{k}|^{2}\Big)^\frac{1}{2}.$
Fig.\ \ref{pic2} reports the numerical solutions and corresponding function $e(t)$ for various fractional order $\alpha$ 
obtained by L1 method with initial values $u_{1}^{0}$ and  $u_{2}^{0}$. 
It clearly shows that the decay rate of $e(t)$ depends directly on the fractional order parameter $\alpha$.
The greater the order $\alpha$, the faster the difference function $e(t)$ contracts. 
But all the contractivity rates remain to be algebraic, rather than the exponential decay rate in the case of  integer-order ODEs 
($\alpha=1$). 

In order to further analyze the quantitative behavior of the decay rate of $e(t)$, we introduce the index:
 \begin{equation}\label{eqn7}
\begin{split}
p_{\alpha}(t) =\frac{\ln\left(c_{\alpha}\|U(0)-V(0)\|\right)-\ln\left(\|U(t)-V(t)\|\right)}{\ln(t)},~~~t >1
\end{split}
\end{equation}
from the improved contractivity rate estimation (\ref{eqn6b}).
Obviously, the index $p_{\alpha}(t) \to\frac{-\ln\left(\|U(t)-V(t)\|\right)}{\ln(t)}$ as $t\to \infty$ and is independent of the initial value $c_{\alpha} \|U(0)-V(0)\|$.
In the numerical simulations, we just take $\|U(1)-V(1)\|= c_{\alpha}\|U(0)-V(0)\|$.

\begin{table}\label{table21}
\caption{The observed index functions $p_\alpha$ by L1 method for Example 4.2 with $h=0.2$} 
\begin{center}
\begin{tabular}{lll lll ll }
\hline $t$  & $\alpha=0.3$  & $\alpha=0.6$ & $\alpha=0.9$ & $\alpha=0.99$\\
\hline
$20   $ & 0.3286  & 0.6771 & 1.0769  & 1.4751\\
$40   $ & 0.3233  &0.6641  & 1.0461  & 1.3866\\
$60  $ & 0.3209  &0.6582  & 1.0324  & 1.3481\\
$80  $& 0.3195  & 0.6546 & 1.0240  & 1.3249\\
$100 $& 0.3185  & 0.6521 & 1.0182  & 1.3089\\
\hline
\end{tabular}
\end{center}
\end{table}

\begin{table}\label{table22}
\caption{The observed index functions $p_\alpha$ by QIA for Example 4.2 with $h=0.2$} 
\begin{center}
\begin{tabular}{lll lll ll }
\hline $t$  & $\alpha=0.3$  & $\alpha=0.6$ & $\alpha=0.9$ & $\alpha=0.99$\\
\hline
$20   $ & 0.3451  & 0.6768 & 0.8355  & 1.7972\\
$40   $ & 0.3375  &0.6639  & 0.8498  & 1.6480\\
$60  $  & 0.3341  &0.6581  &0.8555  & 1.5835\\
$80  $  & 0.3320  &0.6546 & 0.8587  & 1.5449\\
$100 $ & 0.3305  &0.6521 & 0.8609  & 1.5182\\
\hline
\end{tabular}
\end{center}
\end{table}

The observed index $p_\alpha$  for L1 method is presented in Tab.\ \ref{table21} and for QIA method is given in Tab.\ \ref{table22}.
The results show that the contractivity rate is about $\|U(t)-V(t)\|=O(t^{-\alpha})$ as $t\to +\infty$, 
which is consistent with the continuous estimate in (\ref{eqn6b}) 
and our theoretical prediction for numerical contravtivity rate given in \cref{thm:ImpN}.

The semi-discrete F-ODEs (\ref{eqn4}) is stiff when $t$ is small, and F-BDFs work well for relatively large step size $h=0.2$.
As a comparison, when we make use of the F-ABM method proposed in \cite{Diethelm2} for simulations, 
numerical blowup or oscillation appears even for $h=1e-14$ when $\alpha=0.3$. 
The serious restrictions on the step sizes, especially when $\alpha$ is small,  indicate that the explicit F-ABM method is not suitable for stiff F-ODEs. 
In fact, the linear stability of fractional predictor-corrector methods was studied deeply in \cite{Garrappa1}, 
and it was shown that the stability regions of this type of methods are usually relatively small, not suitable for stiff F-ODEs.

\subsection{Nonlinear F-ODEs} Consider  the nonlinear F-ODEs
\begin{equation}\label{eqn8}
\begin{array}{cc}
\begin{split}
(i)~ \left\{
\begin{aligned}
~^{C}_{~0}D_{t}^{\alpha}x(t)&=-x^{3}-x,\\
x(0)&=x_{0}.\\
\end{aligned}
\right.
\end{split}~~~
\begin{split}
(ii)~ \left\{
\begin{aligned}
~^{C}_{~0}D_{t}^{\alpha}x(t)&=-10xy^{2}-x,\\
~^{C}_{~0}D_{t}^{\alpha}y(t)&=10x^{2}y-y,\\
\end{aligned}
\right.
\end{split}
\end{array}
\end{equation}

By simple calculations, it is easy to check that scalar F-ODEs in (i) satisfy the one-sided Lipschitz condition with $\lambda=-1$, 
and also meet the dissipative condition with $a=0, b=1$. Hence, it is contractive and dissipative. 
The F-ODEs in (ii) satisfy the dissipative condition with $a=0, b=1$, hence are dissipative. 

\begin{figure}[!ht]
\begin{center}
$\begin{array}{cc}
 \includegraphics[scale=0.30]{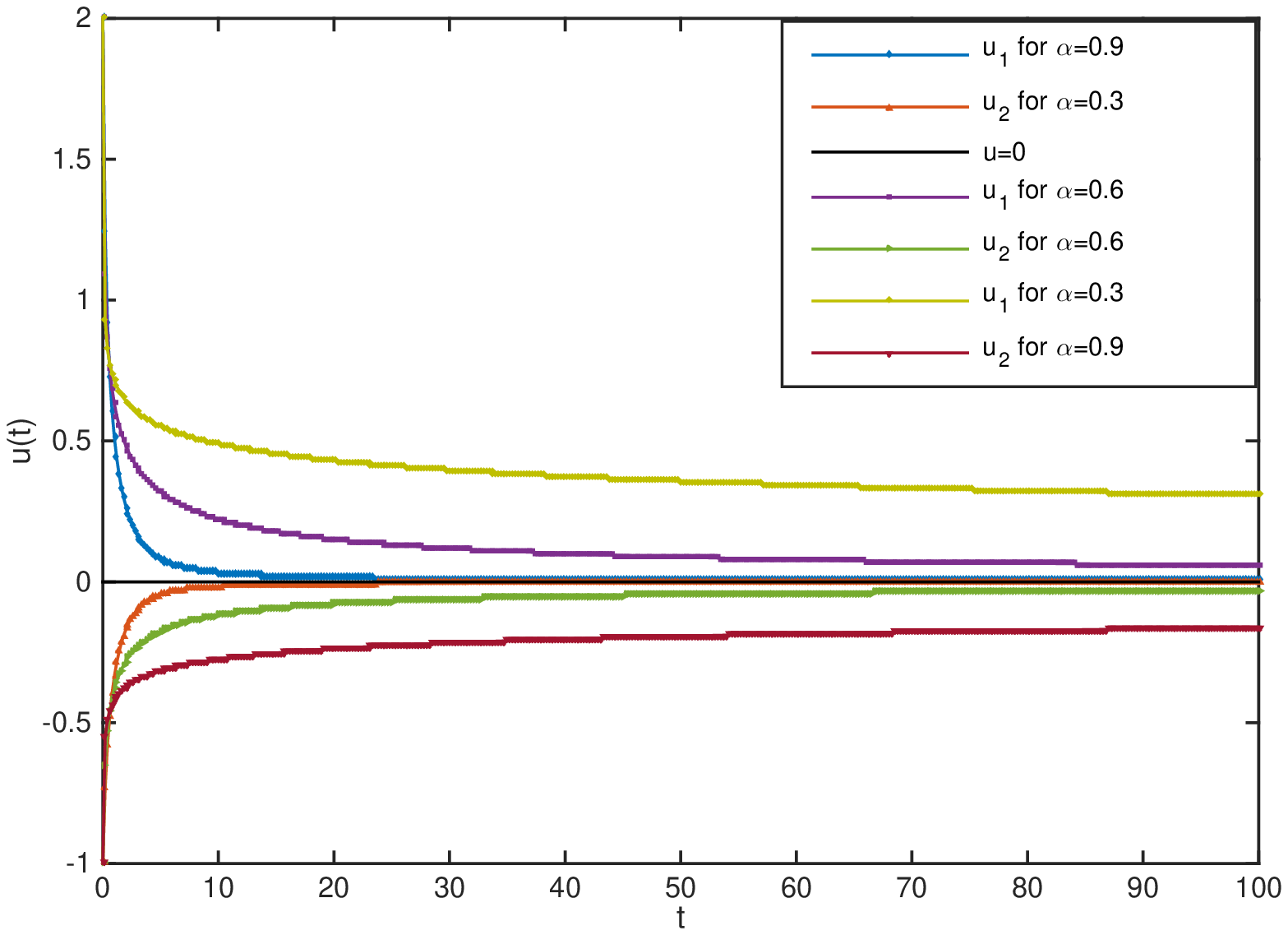}
 \includegraphics[scale=0.35]{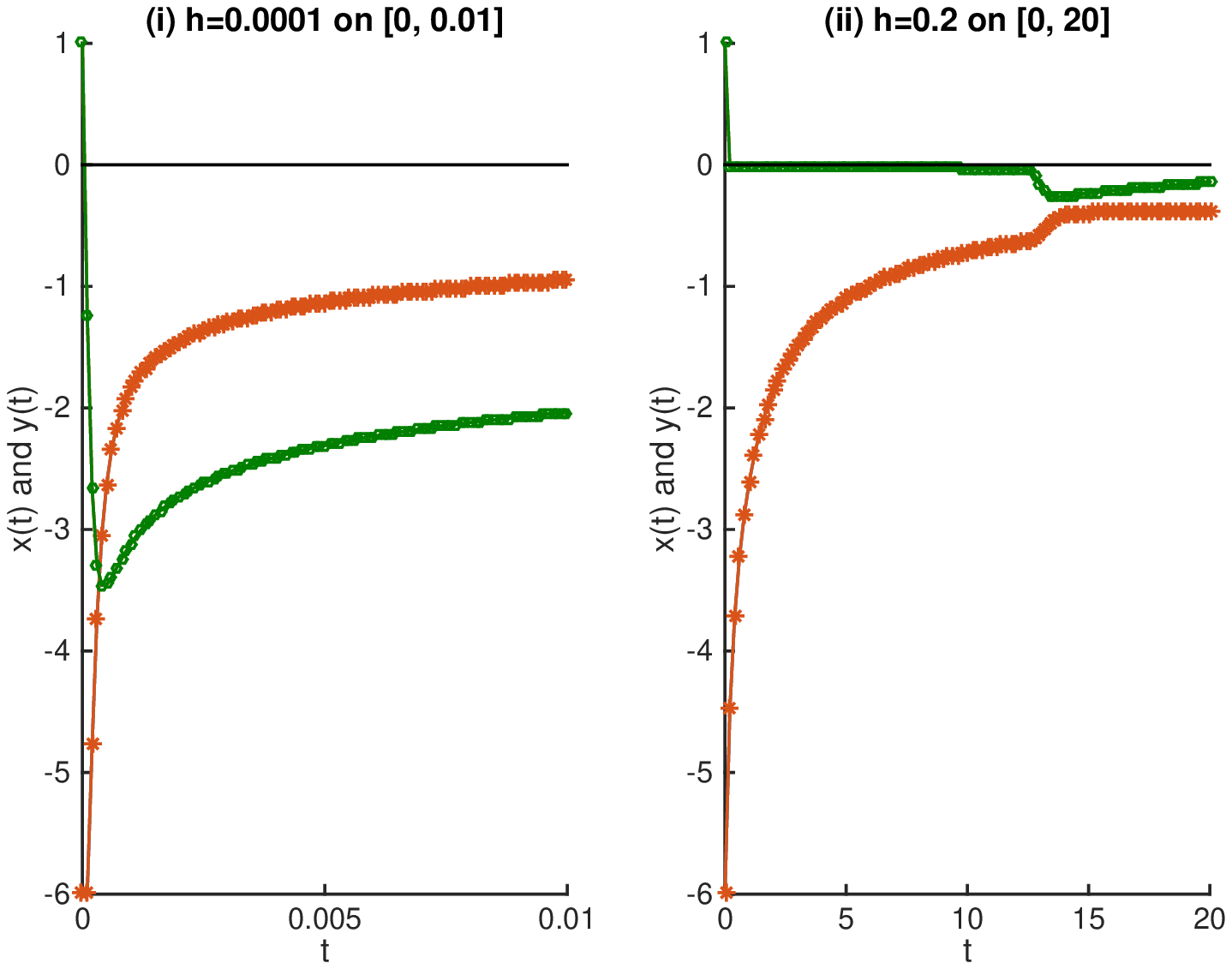}
 \end{array}$
\end{center}
\caption{Left: numerical solutions of (i) at $T=100$ obtained by G-L  with $h=0.2$, initial values $u_{1}=2$, $u_{2}=-1$ and $\alpha=0.3, 0.6, 0.9$; 
Middle: numerical solutions of (ii) at $T=0.01$ obtained by L1  with $h=0.0001$ on $[0, 0.01]$ with $\alpha=0.6$ and initial values $x_{0}=-6, y_{0}=1$;
Right: numerical solutions of (ii) at $T=20$ obtained by L1  with $h=0.2$ on $[0, 20]$ with $\alpha=0.6$ and initial values $x_{0}=-6, y_{0}=1$.}
\label{pic3}
\end{figure}

\begin{table}\label{tabc1}
\caption{Observed index function $p_\alpha$ in Example 4.3 (i) computed by G-L with $h=0.5$ and $T=5000$.} 
\begin{center}
\begin{tabular}{lllllllll }
\hline $t$ & $\alpha=0.3$    & $\alpha=0.6$  & $\alpha=0.9$ & $\alpha=0.99$\\
\hline
$1000   $& 0.2149 &0.5682  &1.0672  & 1.4874\\
$2000   $& 0.2200 &0.5714  &1.0520  & 1.4739\\
$3000  $&  0.2228 &0.5729  &1.0441  & 1.4101\\
$4000  $&  0.2247 &0.5739  & 1.0390 & 1.3906\\
$5000  $&  0.2262 &0.5746  &1.0352  & 1.3755 \\
\hline
\end{tabular}
\end{center}
\end{table}

\begin{table}\label{tabc2} 
\caption{Observed $p_\alpha$ in Example 4.3 (i)  computed by second order F-BDFs with $h=0.5$ and $T=5000$.} 
\begin{center}
\begin{tabular}{lllllllll }
\hline $t$ & $\alpha=0.3$    & $\alpha=0.6$  & $\alpha=0.9$ & $\alpha=0.99$\\
\hline
$1000   $& 0.2333 &0.6036  &1.1183   & 1.5437\\
$2000   $& 0.2367 &0.6035  &1.0984  & 1.4884\\
$3000  $&  0.2387 &0.6035  &1.0882  & 1.4601\\
$4000  $&  0.2407 &0.6034  &1.0817  & 1.4397\\
$5000  $&  0.2412 &0.6034  &1.0767  & 1.4250 \\
\hline
\end{tabular}
\end{center}
\end{table}

\begin{table}\label{tabc3} 
\caption{The observed index function $q_\alpha$ in Example 4.3 (ii) computed by L1 with $h=0.5$ and $T=5000$.}
\begin{center}
\begin{tabular}{lllllllll }
\hline $t$ & $\alpha=0.3$    & $\alpha=0.6$  & $\alpha=0.9$ & $\alpha=0.99$\\
\hline
$1000   $& 0.2662 &0.6038  &1.1094   &1.5342\\
$2000   $& 0.2678 &0.6037  &1.1090  & 1.4830\\
$3000  $&  0.2639 &0.6036  &1.1080  & 1.4552\\
$4000  $&  0.2613 &0.6035  &1.1074  & 1.4362\\
$5000  $&  0.2596 &0.6035  &1.1069  & 1.4246 \\
\hline
\end{tabular}
\end{center}
\end{table}

From Fig. \ref{pic3}, we find that the sign of the numerical solution in scalar F-ODEs (i) 
remains unchanged, as shown in \cref{lem:NumPP}. 
The order $\alpha$ significantly affects the contractivity rate and dissipativity rate, 
and all the solutions decay to zero at a slow rate and keep a long tail.  
As in Example \ref{sec:test2}, we can compute the index $p_{\alpha}$ defined in (\ref{eqn7})  to quantitatively characterize the contractivity rate.
Tab.\,\ref{tabc1} and Tab.\,\ref{tabc2} show that the contractivity rate depends directly 
on the order parameter $\alpha$ and is algebraic,
and almost equal to the rate $\alpha$, which is consistent with the results presented in \cref{thm:ImpN}. 
Note that when $\alpha=0.3$, the index $p_{\alpha}$ is slightly smaller than expected because it takes a long time to get to 
the equilibrium.
For long time simulations, some fast algorithm for Caputo derivatives \cite{JWZhang}, should be very helpful.

From Fig.\,\ref{pic3}, we see that 
the sign of the numerical solutions in the vector F-ODEs (ii) is no longer unchanged. It also exhibits an initial layer and thus has a stiff feature.   
We now introduce an index to quantitatively characterize the dissipativtity rate:
 \begin{equation}\label{eqn9}
\begin{split}
q_{\alpha}(t) =\frac{\ln\left(c_{\alpha}\|u(0)\|\right)-\ln\left(\|u(t)\|\right)}{\ln(t)},~~~t >1.
\end{split}
\end{equation}
The index $q_{\alpha}(t) \to\frac{-\ln\left(\|u(t)\|\right)}{\ln(t)}$ as $t\to +\infty$ and is also independent of the initial values $c_{\alpha} \|u(0)\|$.
In the numerical simulations, we just take $\|u(1)\|= c_{\alpha}\|u(0)\|$.
Tab.\,\ref{tabc3} shows that the dissipativity rate depends directly on the order $\alpha$ and is algebraic,
with the rate nearly equal to $\alpha$.

\section{Concluding remarks}\label{sec:dis}
We have presented some sufficient conditions to ensure 
the numerical contractivity and dissipativity of F-BDFs for nonlinear F-ODEs.
F-BDFs, including four popular schemes, are shown to be dissipative and contractive, 
and can preserve the exact contractivity and dissipativity rates of the solutions to the continuous equations.
To the best of our knowledge, this is the first work on the numerical asymptotic behavior of the solutions to 
nonlinear F-ODEs.
There are still a lot to be done in order to better understand efficient numerical methods for nonlinear F-ODEs 
without the classical Lipschitz conditions.
For instance, stable numerical methods for strongly stiff F-ODEs and their rigorous long-time convergence analysis 
are very important. 
 
We note that for Riemann-Liouville F-ODEs, the numerical dissipativity and contractivity of the F-BDFs that 
we have studied can be developed directly. 
But the decay rate for Riemann-Liouville F-ODEs is slightly changed.
For the multi-order fractional systems with $\alpha=(\alpha_{1}$, $\alpha_{2},...,\alpha_{n})^{T}$, where 
$\alpha_{i}\in(0,1)$ for $i=1,2,...,n$,
their numerical dissipativity and contractivity of F-BDFs can be established in a similar manner.

\appendix
\section{Asymptotical decay rate of Volterra difference equation}
\label{sec:appendixA}
We introduce some related concepts and results in \cite{Applelby}.  
Let $r>0$ be finite. A real sequence $\gamma=\{\gamma_{n}\}_{n\ge0}$ is in $W(r)$ if  $\gamma_{n}> 0$ and
\begin{equation*} 
\begin{split}
\lim_{n\to \infty}\frac{\gamma_{n-1}}{\gamma_{n}}=\frac{1}{r}, ~\tilde{\gamma}(r)=\sum_{i=0}^{\infty}\gamma_{i}r^{-i}<\infty
\quad 
\mbox{and} ~~
\lim_{m\to \infty}\left(\lim_{n\to \infty}\sup \frac{1}{\gamma_{n}} \sum_{i=m}^{n-m}\gamma_{n-i}\gamma_{i} \right)=0.
\end{split}
\end{equation*}

Note that if  $\gamma\in W(r)$ and $r\le 1$, then $\gamma_{n}\to 0$ as $n\to \infty$. The sequence $\gamma_{n}=1/(n+1)^{1+\alpha}\in W(1)$ while 
 $\gamma_{n}=1/(n+1)^{\alpha}$ is not in $W(1)$ for $0<\alpha<1$.
For a given sequence $\gamma=\{\gamma_{n}\}_{n\ge0}$ in $W(r)$ and $x=\{x_{n}\}_{n\ge0}$, we define $L_{\gamma}(x)=\lim\limits_{n\to \infty} \frac{x_{n}}{\gamma_{n}}$ if the limit exists.
We now recall the main results of \cite{Applelby} in the scalar case.

\begin{lemma}[\cite{Applelby}]\label{lem:app}
Consider the Volterra difference equation
\begin{equation} \label{eqa1}
\begin{split}
z_{n+1}=h_{n}+\sum_{i=0}^{n}H_{n, i}z_{i}, ~~n\ge 1.
\end{split}
\end{equation}
Assume that

 \emph{(i)}~ $K:=\lim\limits_{m\to \infty}\sup\left( \lim\limits_{n\to \infty}\sup \sum\limits_{j=0}^{m}|H_{n, n-j}|\right)$ is finite with $K<1$;

\emph{(ii)}~ $H_{n,m}\to H_{\infty, m}$ and  $h_{n}\to h_{\infty}$ as $n\to \infty$ with $\sum\limits_{m=0}^{\infty}|H_{\infty, m}|<\infty$;

\emph{(iii)}~ $\lim\limits_{m\to \infty}\sup\left( \lim\limits_{n\to \infty}\sup \sum\limits_{j=m}^{n-m}|H_{n, j}|\right)=0$.

Then the limit $\lim\limits_{n\to \infty}z_{n}$ exists and satisfies 
$$\lim\limits_{n\to \infty}z_{n}=(1-V)^{-1}\Big(h_{\infty}+\sum\limits_{j=0}^{\infty}H_{\infty, j}z_{j}\Big)
\quad 
\mbox{with} 
~~V:=\lim\limits_{m\to \infty}\Big( \lim\limits_{n\to \infty} \sum\limits_{j=0}^{m}H_{n, n-j}\Big).
$$

Furthermore, when (\ref{eqa1}) is  a convolution equation with $H_{n, j}=H^{\sharp}_{n-j}$ and $\sum\limits_{j=0}^{\infty}|H^{\sharp}_{j}|<1$,
it is seen that $K=\sum\limits_{j=0}^{\infty}|H^{\sharp}_{j}|<1$, $V=\sum\limits_{j=0}^{\infty}H^{\sharp}_{j}$ and $H_{\infty}=0$. Therefore, 
 the limit $\lim\limits_{n\to \infty}z_{n}$ exists and satisfies 
 \begin{equation} \label{eqa2}
\begin{split}
\lim\limits_{n\to \infty}z_{n}=\left(1-\sum\limits_{j=0}^{\infty}H^{\sharp}_{j}\right)^{-1}h_{\infty}.
\end{split}
\end{equation}
\end{lemma}

\section*{Acknowledgements}
The authors are grateful to Professor Zhi Zhou (Hong Kong Polytechnic University) for 
very helpful discussions about some technical issues in this work.

\end{document}

%% file: ex_shared-2.tex
\usepackage{lipsum}
\usepackage{amsfonts}
\usepackage{graphicx}
\usepackage{epstopdf}
\usepackage{algorithmic}
\ifpdf
  \DeclareGraphicsExtensions{.eps,.pdf,.png,.jpg}
\else
  \DeclareGraphicsExtensions{.eps}
\fi

\topmargin by-\headsep  \textheight   230mm   
  \def\nb{\nonumber}
\def\to{\rightarrow}

\newcommand{\q}{\quad}   \newcommand{\qq}{\qquad}

 \def\bs{\bigskip} 

\def\beqn{\begin{eqnarray}}  \def\eqn{\end{eqnarray}}
\def\beqnx{\begin{eqnarray*}} \def\eqnx{\end{eqnarray*}}
\def\bn{\begin{enumerate}} \def\en{\end{enumerate}}

\newcommand{\TheTitle}{Long time behavior of BDFs for F-ODEs} 
\newcommand{\TheAuthors}{Dongling Wang, Aiguo Xiao and Jun Zou}

\headers{\TheTitle}{\TheAuthors}

\title{Long-time behavior of numerical solutions \\ to nonlinear fractional ODEs
}

\author{
  Dongling Wang\thanks{Department of Mathematics and Center for Nonlinear Studies,
Northwest University, Xi'an, Shaanxi 710075, China
    (\email{wdymath@nwu.edu.cn}).
    The work of this author was partially supported by National Natural Science Foundation of China (Grant No.11501447). 
}
  \and
  Aiguo Xiao\thanks{Hunan Key Laboratory for Computation and Simulation in Science and Engineering,
Xiangtan University, Xiangtan,  Hunan 411105, China. The work of this author was partially supported by
National Natural Science Foundation of China (Grant No.11671343). (\email{xag@xtu.edu.cn}).}
  \and
  Jun Zou\thanks{Deptartment of Mathematics, The Chinese University of Hong Kong Shatin, N.T., Hong Kong.
  The work of this author was substantially supported by Hong Kong RGC General Research Fund (projects 405513 and 14306814). 
    (\email{zou@math.cuhk.edu.hk}).}
  }
  
